\documentclass[onefignum,onetabnum]{siamonline171218}
\usepackage{amsmath,amsfonts,amssymb}
\usepackage{dsfont}
\usepackage{array}
\usepackage{stfloats}
\usepackage{subcaption}
\DeclareMathSizes{10}{10}{5}{5}

\usepackage{lipsum}
\usepackage{amsfonts}
\usepackage{graphicx}
\usepackage{epstopdf}
\usepackage{algorithmic}
\ifpdf
  \DeclareGraphicsExtensions{.eps,.pdf,.png,.jpg}
\else
  \DeclareGraphicsExtensions{.eps}
\fi

\usepackage{enumitem}
\setlist[enumerate]{leftmargin=.5in}
\setlist[itemize]{leftmargin=.5in}

\newsiamremark{remark}{Remark}
\newsiamremark{hypothesis}{Hypothesis}
\crefname{hypothesis}{Hypothesis}{Hypotheses}
\newsiamthm{claim}{Claim}

\headers{Recursive estimation in a Riemannian manifold}{J. Zhou, and S. Said}

\title{Fast, asymptotically efficient, recursive \\estimation in a Riemannian manifold\thanks{Submitted to the editors 15 July 2018.
\funding{This work was funded by Cluster SysNum, part of the Excellence Initiative of Universit\'e de Bordeaux.}}}

\author{Jialun Zhou\thanks{Universit\'e de Franche-Comt\'e, Master Mod\'elisation Statistique (\email{jialun.zhou@edu.univ-fcomte.fr}).}
\and Salem Said\thanks{CNRS, Universit\'e de Bordeaux, (Laboratoire IMS UMR 5218) (\email{salem.said@u-bordeaux.fr}).}}

\usepackage{amsopn}

\ifpdf
\hypersetup{
  pdftitle={Fast, asymptotically efficient, recursive estimation in a Riemannian manifold},
  pdfauthor={J. Zhou, and S. Said}
}
\fi

\begin{document}

\maketitle

\begin{abstract}
Stochastic optimisation in Riemannian manifolds, especially the Riemannian stochastic gradient method, has attracted much recent attention. The present work applies stochastic optimisation to the task of recursive estimation of a statistical parameter which belongs to a Riemannian manifold. Roughly, this task amounts to stochastic minimisation of a statistical divergence function. The following problem is considered\,: how to obtain fast, asymptotically efficient, recursive estimates, using a Riemannian stochastic optimisation algorithm with decreasing step sizes? In solving this problem, several original results are introduced. First, without any convexity assumptions on the divergence function, it is proved that, with an adequate choice of step sizes, the algorithm computes recursive estimates which achieve a fast non-asymptotic rate of convergence. Second, the asymptotic normality of these recursive estimates is proved, by employing a novel linearisation technique. Third, it is proved that, when the Fisher information metric is used to guide the algorithm, these recursive estimates achieve an optimal asymptotic rate of convergence, in the sense that they become asymptotically efficient. These results, while relatively familiar in the Euclidean context, are here formulated and proved for the first time, in the Riemannian context. In addition, they are illustrated with a numerical application to the recursive estimation of elliptically contoured distributions. 
\end{abstract}  


\section{Introduction}\label{sec:intro}
Over the last five years, the data science community has devoted significant attention to stochastic optimisation in Riemannian manifolds. This was impulsed by Bonnabel, who proved the convergence of the Riemannian stochastic gradient method~\cite{bonnabel}. Later on~\cite{sra}, the rate of convergence of this method was studied in detail, under various convexity assumptions on the cost function. More recently, asymptotic efficiency of the averaged Riemannian stochastic gradient method was proved in~\cite{flamm}. Previously, for the specific problem of computing Riemannian means, several results on the convergence and asymptotic normality of Riemannian stochastic optimisation methods had been obtained~\cite{arnaudon}\cite{yang}. 
   
The present work moves in a different direction, focusing on recursive estimation in Riemannian manifolds. While recursive estimation is a special case of stochastic optimisation, it has its own geometric structure, given by the Fisher information metric. Here, several original results will be introduced, which show how this geometric structure can be exploited, to design Riemannian stochastic optimisation algorithms which compute fast, asymptotically efficient, recursive estimates, of a statistical parameter which belongs to a Riemannian manifold. For the first time in the literature, these results extend, from the Euclidean context to the Riemannian context, the classical results of~\cite{nev}\cite{duflo}. 

The mathematical problem, considered in the present work, is formulated in Section \ref{sec:problem}. This involves a parameterised statistical model $P$ of probability distributions $P_{\scriptscriptstyle \theta\,}$, where the statistical parameter $\theta$ belongs to a Riemannian manifold $\Theta$. Given independent observations, with distribution $P_{\scriptscriptstyle \theta^*}$ for some $\theta^* \in \Theta$, the aim is to estimate the unknown parameter $\theta^*$. In principle, this is done by minimising a statistical divergence function $D(\theta)$, which measures the dissimilarity between $P_{\scriptscriptstyle \theta}$ and $P_{\scriptscriptstyle \theta^*\,}$. Taking advantage of the observations, there are two approaches to minimising $D(\theta)$\,: stochastic minimisation, which leads to recursive estimation, and empirical minimisation, which leads to classical techniques, such as maximum-likelihood estimation~\cite{broniatowski2}\cite{broniatowski1}. 

The original results, obtained in the present work, are stated in Section \ref{sec:props}. In particular, these are Propositions \ref{prop:ratel2}, \ref{prop:normality}, and \ref{prop:information}. Overall, these propositions show that recursive estimation, which requires less computational resources than maximum-likelihood estimation, can still achieve the same optimal performance, characterised by asymptotic efficiency~\cite{ibrahas}\cite{vaart}. 

To summarise these propositions, consider a sequence of recursive estimates $\theta_{\scriptscriptstyle n\,}$, computed using a Riemannian stochastic optimisation algorithm with decreasing step sizes ($n$ is the number of observations already processed by the algorithm). Informally, under assumptions which guarantee that $\theta^*$ is an attractive local minimum of $D(\theta)$, and that the algorithm is neither too noisy, nor too unstable, in the neighborhood of $\theta^*$,\\
\indent $\bullet$ Proposition \ref{prop:ratel2} states that, with an adequate choice of step sizes, the $\theta_{\scriptscriptstyle n}$ achieve a fast non-asymptotic rate of convergence to $\theta^*$. Precisely, the expectation of the squared Riemannian distance between $\theta_{\scriptscriptstyle n}$ and $\theta^*$ is $O\left(n^{\scriptscriptstyle -1}\right)$. This is called a fast rate, because it is the best achievable, for any step sizes which are proportional 
to $n^{\scriptscriptstyle -q}$ with $q \in (1/2,1]$~\cite{benv}\cite{duflo}. Here, this rate is obtained without any convexity assumptions, for twice differentiable $D(\theta)$. It would still hold for non-differentiable, but strongly convex, $D(\theta)$~\cite{sra}. \\
\indent $\bullet$  Proposition \ref{prop:normality} states that the distribution of the $\theta_{\scriptscriptstyle n}$ becomes asymptotically normal, centred at $\theta^*$, when $n$ grows increasingly large, and also characterises the corresponding asymptotic covariance matrix. This proposition is proved using a novel linearisation technique, which also plays a central role in~\cite{flamm}. \\
\indent $\bullet$ Proposition \ref{prop:information} states that, if the Riemannian manifold $\Theta$ is equipped with the Fisher information metric of the statistical model $P$, then Riemannian gradient descent with respect to this information metric, when used to minimise $D(\theta)$, computes recursive estimates $\theta_{\scriptscriptstyle n}$ which are asymptotically efficient, achieving the optimal asymptotic rate of convergence, given by the Cram\'er-Rao lower bound.  This is illustrated, with a numerical application to the recursive estimation of elliptically contoured distributions, in Section \ref{sec:mggd}. 
  
\indent The end result of Proposition \ref{prop:information} is asymptotic efficiency, achieved using the Fisher information metric. In~\cite{flamm}, an alternative route to asymptotic efficiency is proposed, using the averaged Riemannian stochastic gradient method. This method does not require any prior knowledge of the Fisher information metric, but has an additional computational cost, which comes from computing on-line Riemannian averages. 

The proofs of Propositions \ref{prop:ratel2}, \ref{prop:normality}, and \ref{prop:information}, are detailed in Section \ref{sec:proofs}, and Appendices \ref{sec:geometric} and \ref{sec:clt}. Necessary background, about the Fisher information metric (in short, this will be called the information metric),  is recalled in Appendix \ref{sec:efficiency}. Before going on, the reader should note that the summation convention of differential geometry is used throughout the following, when working in local coordinates.


%


\vfill
\pagebreak

\section{Problem statement} \label{sec:problem}Let $P = (P,\Theta,X)$ be a statistical model, with parameter space $\Theta$ and sample space $X$. To each $\theta \in \Theta$, the model $P$ associates a probability distribution $P_{\scriptscriptstyle\theta}$ on $X$. Here, $\Theta$ is a $C^r$ Riemannian manifold with $r > 3$, and $X$ is any measurable space. The Riemannian metric of $\Theta$ will be denoted $\langle\cdot,\cdot\rangle$, with its Riemannian distance $d(\cdot,\cdot)$. In general, the metric $\langle\cdot,\cdot\rangle$ is not the information metric of the model $P$.

Let $(\Omega,\mathcal{F},\mathbb{P})$ be a complete probability space, and $(x_{\scriptscriptstyle n}\,;n=1,2,\ldots)$ be i.i.d. random variables on 
$\Omega$, with values in $X$. While the distribution of $x_{\scriptscriptstyle n}$ is unknown, it is assumed to belong to the model $P$. That is, $\mathbb{P}\circ x^{\scriptscriptstyle -1}_{\scriptscriptstyle n} = P_{\scriptscriptstyle \theta^*}$ for some $\theta^* \in \Theta$, to be called the true parameter. 

Consider the following problem\,: how to obtain fast, asymptotically efficient, recursive estimates $\theta_{\scriptscriptstyle n}$ of the true parameter $\theta^*$, based on observations of the random variables $x_{\scriptscriptstyle n}$? The present work proposes to solve this problem through a detailed study of the decreasing-step-size algorithm, which computes
\begin{subequations} \label{subeq:algorithm}
\begin{equation} \label{eq:algorithm}
   \theta_{\scriptscriptstyle n+1} = \mathrm{Exp}_{\scriptscriptstyle \theta_{\scriptscriptstyle n}}\!\left(\gamma_{\scriptscriptstyle n+1}u(\theta_{\scriptscriptstyle n},x_{\scriptscriptstyle n+1})\right) \hspace{1cm} n = 0,1,\ldots
\end{equation}
starting from an initial guess $\theta_{\scriptscriptstyle 0}\,$. 

This algorithm has three ingredients. First, $\mathrm{Exp}$ denotes the Riemannian exponential map of the metric $\langle\cdot,\cdot\rangle$ of $\Theta$~\cite{petersen}. Second, the step sizes $\gamma_{\scriptscriptstyle n}$ are strictly positive, decreasing, and verify the usual conditions for stochastic approximation~\cite{nev}\cite{kushner}
\begin{equation} \label{eq:stepsize}
    \sum\,\gamma_{\scriptscriptstyle n} \,=\, \infty \hspace{1cm}   \sum\,\gamma^{\scriptscriptstyle 2}_{\scriptscriptstyle n} \,<\, \infty
\end{equation}
Third, $u(\theta,x)$ is a continuous vector field on $\Theta$ for each $x \in X$, which generalises the classical concept of score statistic~\cite{ibrahas}\cite{heyde}. It will become clear, from the results given in Section \ref{sec:props}, that the solution of the above-stated problem depends on the choice of each one of these three ingredients. 

A priori knowledge about the model $P$ is injected into Algorithm (\ref{eq:algorithm}) using a divergence function $D(\theta) = D(P_{\scriptscriptstyle \theta^*},P_{\scriptscriptstyle \theta})$. As defined in~\cite{amari}, this is a positive function, equal to zero if and only if $P_{\scriptscriptstyle \theta} = P_{\scriptscriptstyle \theta^*\,}$, and with positive definite Hessian at $\theta = \theta^*$. Since one expects that minimising $D(\theta)$ will lead to estimating $\theta^*$, it is natural to require that 
\begin{equation} \label{eq:gradient}
   E_{\scriptscriptstyle\theta^*\,}u(\theta,x) \,=\, - \nabla D(\theta)
\end{equation}
In other words, that $u(\theta,x)$ is an unbiased estimator of minus the Riemannian gradient of $D(\theta)$. With $u(\theta,x)$ given by (\ref{eq:gradient}), Algorithm (\ref{eq:algorithm}) is a Riemannian stochastic gradient descent, of the form considered in~\cite{bonnabel}\cite{sra}\cite{flamm}. However, as explained in Remark \ref{rk:gradient}, (\ref{eq:gradient}) may be replaced by the weaker condition (\ref{eq:weakgrad1}), without affecting the results in Section \ref{sec:props}. In this sense, Algorithm (\ref{eq:algorithm}) is more general than Riemannian stochastic gradient descent.  
\end{subequations}  

In practice, a suitable choice of $D(\theta)$ is often the Kullback-Leibler divergence~\cite{shiryayev}, 
\begin{subequations} \label{subeq:kl}
\begin{equation} \label{eq:kl}
   D(\theta) \,=\, -\,E_{\scriptscriptstyle\theta^*}\log L(\theta) \hspace{1cm} L(\theta) \,=\, \frac{dP_{\scriptscriptstyle \theta}}{dP_{\scriptscriptstyle \theta^*}}
\end{equation}
where $P_{\scriptscriptstyle \theta}$ is absolutely continuous with respect to $P_{\scriptscriptstyle \theta^*}$ with Radon-Nikodym derivative $L(\theta)$.
Indeed, if $D(\theta)$ is chosen to be the Kullback-Leibler divergence, then (\ref{eq:gradient}) is satisfied by
\begin{equation} \label{eq:score}
  u(\theta,x) = \nabla \log L(\theta)
\end{equation}
which, in many practical situations, can be evaluated directly, without any knowledge of $\theta^*\,$.
\end{subequations}

\section{Main results} \label{sec:props}The motivation of the following Propositions \ref{prop:as} to \ref{prop:information} is to provide general conditions, which guarantee that Algorithm (\ref{eq:algorithm}) computes fast, asymptotically efficient, recursive estimates $\theta_{\scriptscriptstyle n}$ of the true parameter $\theta^*$. In the statement of these propositions, it is implicitly assumed that conditions (\ref{eq:stepsize}) and (\ref{eq:gradient}) are verified. Moreover, the following assumptions are considered. \\[0.1cm]
\indent \textbf{(d1)} the divergence function $D(\theta)$ has an isolated stationary point at $\theta = \theta^*$, and Lipschitz gradient in a neighborhood of this point. 

\textbf{(d2)} this stationary point is moreover attractive\,: $D(\theta)$ is twice differentiable at $\theta = \theta^*$, with positive definite Hessian at this point.

\textbf{(u1)} in a neighborhood of $\theta = \theta^*$, the function $V(\theta) = E_{\scriptscriptstyle \theta^*}\Vert u(\theta,x)\Vert^2$ is uniformly bounded.

\textbf{(u2)} in a neighborhood of $\theta = \theta^*$, the function $R(\theta) = E_{\scriptscriptstyle \theta^*}\Vert u(\theta,x)\Vert^4$ is uniformly bounded. \\[0.1cm]
For Assumption (d1), the definition of a Lipschitz vector field on a Riemannian manifold may be found in~\cite{meunier}. For Assumptions (u1) and (u2), $\Vert\cdot\Vert$ denotes the Riemannian norm.\\[0.1cm] 
\indent Let $\Theta^*$ be a neighborhood of $\theta^*$ which verifies (d1), (u1), and (u2). Without loss of generality, it is assumed that $\Theta^*$ is compact and convex (see the definition of convexity in~\cite{petersen}\cite{udriste}). Then, $\Theta^*$ admits a system of normal coordinates $(\theta^{\scriptscriptstyle\,\alpha}\,;\alpha = 1\,,\ldots,\,d\,)$ with origin at $\theta^*$. With respect to these coordinates, denote the components of $u(\theta^*,x)$ by $u^{\scriptscriptstyle \alpha}(\theta^*)$ and let $\Sigma^* = (\Sigma^*_{\scriptscriptstyle\alpha\beta})$,\begin{subequations}
\begin{equation} \label{eq:cov}
  \Sigma^*_{\scriptscriptstyle\alpha\beta} \,=\, E_{\scriptscriptstyle \theta^*}\! \left[u^{\scriptscriptstyle\alpha}(\theta^*)\,u^{\scriptscriptstyle\beta}(\theta^*)\right]
\end{equation}
When (d2) is verified, denote the components of the Hessian of $D(\theta)$ at $\theta = \theta^*$ by $H = \left(H_{\scriptscriptstyle \alpha \beta}\right)$, 
\begin{equation} \label{eq:hess}
H_{\scriptscriptstyle \alpha \beta} \,=\, \left.\frac{\partial^{\scriptscriptstyle\, 2}\!\,D}{\mathstrut\partial\theta^{\scriptscriptstyle \alpha}\partial\theta^{\scriptscriptstyle \beta}}\right|_{\scriptscriptstyle \theta^{\scriptscriptstyle \alpha} = 0} 
\end{equation}
Then, the matrix $H = \left(H_{\scriptscriptstyle \alpha \beta}\right)$ is positive definite~\cite{absil}. Denote by $\lambda > 0$ its smallest eigenvalue.
\end{subequations}

Propositions \ref{prop:as} to \ref{prop:information} require the condition that the recursive estimates $\theta_{\scriptscriptstyle n}$ are stable, which means that all the $\theta_n$ lie in $\Theta^*$, almost surely. The need for this condition is discussed in Remark \ref{rk:stable}. Note that, if $\theta_{\scriptscriptstyle n}$ lies in $\Theta^*$, then $\theta_{\scriptscriptstyle n}$ is determined by its normal coordinates $\theta^{\scriptscriptstyle\, \alpha}_{\scriptscriptstyle n\,}$. 
\begin{proposition}[consistency]\label{prop:as}assume (d1) and (u1) are verified, and the recursive estimates $\theta_{\scriptscriptstyle n}$ are stable. Then, $\lim\theta_{\scriptscriptstyle n} = \theta^*$ almost surely. 
\end{proposition}

\begin{proposition}[mean-square rate]\label{prop:ratel2}assume (d1), (d2) and (u1) are verified, the recursive estimates $\theta_{\scriptscriptstyle n}$ are stable, and $\gamma_{\scriptscriptstyle n} = \frac{a}{n}$ where $2\lambda a > 1$. Then
\begin{equation} \label{eq:ratel2}
     \mathbb{E}\,d^{\scriptscriptstyle\, 2}(\theta_{\scriptscriptstyle n\,},\theta^*) \,=\, O\left(n^{\scriptscriptstyle -1}\right)
\end{equation}
\end{proposition}
\begin{proposition}[almost-sure rate]\label{prop:rateas}assume the conditions of Proposition \ref{prop:ratel2} are verified. Then,
\begin{equation} \label{eq:rateas}
d^{\scriptscriptstyle\, 2}(\theta_{\scriptscriptstyle n\,},\theta^*) \,=\, o(n^{\scriptscriptstyle -p}) \text{ for } p \in (0,1) \hspace{1cm} \text{almost surely}
\end{equation}
\end{proposition}
\begin{proposition}[asymptotic normality]\label{prop:normality}assume the conditions of Proposition \ref{prop:ratel2}, as well as (u2), are verified. Then, the distribution of the re-scaled coordinates $(n^{\scriptscriptstyle 1/2}\theta^{\scriptscriptstyle\,\alpha}_{\scriptscriptstyle n})$ converges to a centred $d$-variate normal distribution, with covariance matrix $\Sigma$ given by Lyapunov's equation
\begin{equation} \label{eq:lyapunov}
 A\,\Sigma\,+\Sigma\,A  \,=\, -a^{\scriptscriptstyle 2}\,\Sigma^*
\end{equation}
where $A = \left(A_{\scriptscriptstyle \alpha\beta}\right)$ with $A_{\scriptscriptstyle \alpha\beta} = \frac{1}{2}\delta_{\scriptscriptstyle \alpha\beta} - aH_{\scriptscriptstyle \alpha \beta}$ (here, $\delta$ denotes Kronecker's delta).
\end{proposition}
\begin{proposition}[asymptotic efficiency]\label{prop:information}assume the Riemannian metric $\langle\cdot,\cdot\rangle$ of $\Theta$ coincides with the information metric of the model $P$, and let $D(\theta)$ be the Kullback-Leibler divergence (\ref{eq:kl}). Further, assume (d1), (d2), (u1) and (u2) are verified, the recursive estimates $\theta_{\scriptscriptstyle n}$ are stable, and  $\gamma_{\scriptscriptstyle n} = \frac{a}{n}$ where $2a > 1$. Then,\\[0.1cm]
\begin{subequations} \label{subeq:information}
\noindent \textbf{\emph{(i)}} the rates of convergence (\ref{eq:ratel2}) and (\ref{eq:rateas}) hold true. \\[0.1cm]
\noindent \textbf{\emph{(ii)}} if $a = 1$, the distribution of the re-scaled coordinates $(n^{\scriptscriptstyle 1/2}\theta^{\scriptscriptstyle\,\alpha}_{\scriptscriptstyle n})$ converges to a centred $d$-variate normal distribution, with covariance matrix $\Sigma^*$.\\[0.1cm]
\noindent \textbf{\emph{(iii)}} if $a = 1$, and $u(\theta,x)$ is given by (\ref{eq:score}), then $\Sigma^*$ is the identity matrix, and the recursive estimates $\theta_{\scriptscriptstyle n}$ are asymptotically efficient. \\[0.1cm]
\noindent \textbf{\emph{(iv)}} the following rates of convergence also hold
\begin{eqnarray} 
\label{eq:information1} \mathbb{E}\,D(\theta_{\scriptscriptstyle n})\,=\, O\left(n^{\scriptscriptstyle -1}\right) \hspace{4.9cm} \\[0.1cm]
\label{eq:information2} D(\theta_{\scriptscriptstyle n\,}) \,=\, o(n^{\scriptscriptstyle -p}) \text{ for } p \in (0,1)  \hspace{1cm} \text{almost surely}
\end{eqnarray}
\end{subequations}
\end{proposition}
The following remarks are concerned with the scope of Assumptions (d1), (d2), (u1), and (u2), and with the applicability of Propositions \ref{prop:as} to \ref{prop:information}. 
\begin{remark}\label{rk:metric}(d2), (u1) and (u2) do not depend on the Riemannian metric $\langle\cdot,\cdot\rangle$ of $\Theta$. Precisely, if they are verified for one Riemannian metric on $\Theta$, then they are verified for any Riemannian metric on $\Theta$. Moreover, if the function $D(\theta)$ is $C^{\scriptscriptstyle2}$, then the same is true for (d1). In this case, Propositions \ref{prop:as} to \ref{prop:information} apply for any Riemannian metric on $\Theta$, so that the choice of the metric $\langle\cdot,\cdot\rangle$ is a purely practical matter, to be decided according to applications. 
\end{remark}
\begin{remark}\label{rk:gradient}the conclusion of Proposition \ref{prop:as} continues to hold, if (\ref{eq:gradient}) is replaced by 
\begin{equation} \label{eq:weakgrad1}
E_{\scriptscriptstyle \theta^*}\langle u(\theta,x),\!\nabla D(\theta)\rangle < 0 \text{ for } \theta \neq \theta^*
\end{equation}
Then, it is even possible to preserve Propositions \ref{prop:ratel2}, \ref{prop:rateas}, and \ref{prop:normality}, provided (d2) is replaced by the assumption that the mean vector field, $X(\theta) = E_{\scriptscriptstyle \theta^*\,} u(\theta,x)$, has an attractive stationary point at $\theta = \theta^*$. This generalisation of Propositions \ref{prop:as} to \ref{prop:normality} can be achieved following essentially the same approach as laid out in Section \ref{sec:proofs}. However, in the present work, it will not be carried out in detail.
\end{remark}
\begin{remark}\label{rk:stable}the condition that the recursive estimates $\theta_{\scriptscriptstyle n}$ are stable is standard in all prior work on stochastic optimisation in manifolds~\cite{bonnabel}\cite{sra}\cite{flamm}. In practice, this condition can be enforced through replacing Algorithm (\ref{eq:algorithm}) by a so-called projected or truncated algorithm. This is identical to (\ref{eq:algorithm}), except that $\theta_{\scriptscriptstyle n}$ is projected back onto the neighborhood $\Theta^*$ of $\theta^*$, whenever it falls outside of this neighborhood~\cite{nev}\cite{kushner}. On the other hand, if the $\theta_{\scriptscriptstyle n}$ are not required to be stable, but (d1) and (u1) are replaced by global assumptions,
\\[0.1cm]
\indent \textbf{(d1')} $D(\theta)$ has compact level sets and globally Lipschitz gradient. \\[0.1cm]
\indent \textbf{(u1')} $V(\theta) \leq C\,(1+D(\theta))$ for some constant $C$ and for all $\theta \in \Theta$. \\[0.1cm]
then, applying the same arguments as in the proof of Proposition \ref{prop:as}, it follows that the $\theta_{\scriptscriptstyle n}$ converge to the set of stationary points of $D(\theta)$, almost surely. 
\end{remark}
\begin{remark}\label{rk:chi2} from (ii) and (iii) of Proposition \ref{prop:information}, it follows that the distribution of $n\,d^{\scriptscriptstyle\,2}(\theta_{\scriptscriptstyle n\,},\theta^*)$ converges to a $\chi^{\scriptscriptstyle 2}$-distribution with $d$ degrees of freedom. This provides a practical means of confirming the asymptotic efficiency of the recursive estimates $\theta_{\scriptscriptstyle n\,}$.
\end{remark}
   
\section{Application\,: estimation of ECD} \label{sec:mggd}Here, the conclusion of Proposition \ref{prop:information} is illustrated, by
applying Algorithm (\ref{eq:algorithm}) to the estimation of elliptically contoured distributions (ECD) \cite{kotz}\cite{sraell}. Precisely, in the notation of Section \ref{sec:problem}, let $\Theta = \mathcal{P}_{\scriptscriptstyle m}$ the space of $m \times m$ positive definite matrices, and $X = \mathbb{R}^{\scriptscriptstyle m\,}$. Moreover, let each $P_{\scriptscriptstyle \theta}$ have probability density function
\begin{equation} \label{eq:ecd}
  p(x|\theta) \,\propto\, \exp\left[ h\left(x^{\scriptscriptstyle{\dagger}}\theta^{\scriptscriptstyle-1}x\right) - \frac{1}{2}\log\det(\theta)\right] \hspace{1cm} \theta \in \mathcal{P}_{\scriptscriptstyle m}\,,x \in \mathbb{R}^{\scriptscriptstyle m}
\end{equation}
where $h:\mathbb{R}\rightarrow \mathbb{R}$ is fixed, has negative values, and is decreasing, and $^{\scriptscriptstyle\dagger}$ denotes the transpose. Then, $P_{\scriptscriptstyle \theta}$ is called an ECD with scatter matrix $\theta$. To begin, let $(x_{\scriptscriptstyle n}\,;n=1,2,\ldots)$ be i.i.d. random vectors in $\mathbb{R}^{\scriptscriptstyle m\,}$, with distribution $P_{\scriptscriptstyle\theta^*}$ given by (\ref{eq:ecd}), and consider the problem of estimating the true scatter matrix $\theta^*$. The standard approach to this problem is based on maximum-likelihood estimation~\cite{pascal}\cite{sraell}. An original approach, based on recursive estimation, is now introduced using Algorithm (\ref{eq:algorithm}).

As in Proposition \ref{prop:information}, the parameter space $\mathcal{P}_{\scriptscriptstyle m}$ will be equipped with the information metric of the statistical model $P$ just described. In~\cite{berkane}, it is proved that this information metric is an affine-invariant metric on $\mathcal{P}_{\scriptscriptstyle m\,}$. In other words, it is of the general form~\cite{cyrus} 
\begin{subequations} \label{subeq:infometric}
\begin{equation} \label{eq:affinv}
   \langle u,u\rangle_{\scriptscriptstyle \theta} \,=\,  I_{\scriptscriptstyle 1}\,\mathrm{tr}\left(\theta^{\scriptscriptstyle -1}u\right)^{\scriptscriptstyle 2}\,+\,I_{\scriptscriptstyle 2}\,\mathrm{tr}^{\scriptscriptstyle 2}\left(\theta^{\scriptscriptstyle -1}u\right) \hspace{1cm} u \in T_{\scriptscriptstyle \theta}\mathcal{P}_{\scriptscriptstyle m}
\end{equation}
parameterised by constants $I_{\scriptscriptstyle 1} > 0$ and $I_{\scriptscriptstyle 2} \geq 0$, where $\mathrm{tr}$ denotes the trace and $\mathrm{tr}^{\scriptscriptstyle 2}$ the squared trace. Precisely~\cite{berkane}, for the information metric of the model $P$, 
\begin{equation} \label{eq:infocoeff}
  I_{\scriptscriptstyle 1} = \frac{\varphi}{2m^{\scriptscriptstyle 2}(m+2)} \hspace{1cm} 
  I_{\scriptscriptstyle 2} = \frac{\varphi}{m^{\scriptscriptstyle 2}} - \frac{1}{4} 
\end{equation} 
where $\varphi$ is a further constant, given by the expectation
\begin{equation} \label{eq:varphi}
\varphi \,=\, E_{\scriptscriptstyle e}\left[h^\prime(x^{\scriptscriptstyle{\dagger}}x)\left(x^{\scriptscriptstyle{\dagger}}x\right)\right]^2
\end{equation}
with $e \in \mathcal{P}_{\scriptscriptstyle m}$ the identity matrix, and $h^\prime$ the derivative of $h$. This expression of the information metric can now be used to specify Algorithm (\ref{eq:algorithm}).
\end{subequations}

First, since the information metric is affine-invariant, it is enough to recall that all affine-invariant metrics on $\mathcal{P}_{\scriptscriptstyle m}$ have the same Riemannian exponential map~\cite{pennec}\cite{sraell},
\begin{subequations} \label{subeq:algoecd}
  \begin{equation} \label{eq:exp}
   \mathrm{Exp}_{\scriptscriptstyle \theta}(u) \,=\, \theta\exp\left(\theta^{\scriptscriptstyle -1}u\right)
 \end{equation}  
where $\exp$ denotes the matrix exponential. Second, as in (ii) of Proposition \ref{prop:information}, choose the sequence of step sizes
\begin{equation} \label{eq:assstep}
\gamma_{\scriptscriptstyle n} = \frac{1}{n}
\end{equation}
Third, as in (iii) of Proposition \ref{prop:information}, let $u(\theta,x)$ be the vector field on $\mathcal{P}_{\scriptscriptstyle m}$ given by (\ref{eq:score}), 
\begin{equation} \label{eq:scoreecd}
  u(\theta,x) = \nabla^{\scriptscriptstyle (inf)} \log L(\theta) = \nabla^{\scriptscriptstyle (inf)} \log p(x|\theta)
\end{equation}  
where $\nabla^{\scriptscriptstyle (inf)}$ denotes the gradient with respect to the information metric, and $L(\theta)$ is the likelihood ratio, equal to $p(x|\theta)$ divided by $p(x|\theta^*)$. Now, replacing (\ref{subeq:algoecd}) into (\ref{eq:algorithm}) defines an original algorithm for recursive estimation of the true scatter matrix $\theta^*$.
\end{subequations}

To apply this algorithm in practice, one may evaluate $u(\theta,x)$ via the following steps. Denote $g(\theta,x)$ the gradient of $\log p(x|\theta)$ with respect to the affine-invariant metric of~\cite{pennec}, which corresponds to $I_{\scriptscriptstyle 1} = 1$ and $I_{\scriptscriptstyle 2} = 0$. By direct calculation from (\ref{eq:ecd}), this is given by
 \begin{subequations} \label{subeq:infograd}
\begin{equation} \label{eq:classgrad}
g(\theta,x) \,=\, -\frac{1}{2}\theta - h^\prime\left(x^{\scriptscriptstyle\dagger}\theta^{\scriptscriptstyle -1}x\right)x x^{\scriptscriptstyle\dagger}
\end{equation}
Moreover, introduce the constants $J_{\scriptscriptstyle 1}= I_{\scriptscriptstyle 1}$ and $J_{\scriptscriptstyle 2} = I_{\scriptscriptstyle 1}+ m I_{\scriptscriptstyle 2\,}$. Then, $u(\theta,x)$ can be evaluated,
\begin{equation} \label{eq:uthetax}
  u(\theta,x) \,=\, J^{\scriptscriptstyle -1}_{\scriptscriptstyle 1} \left(g(\theta,x)\right)^{\scriptscriptstyle \perp}\,+\,J^{\scriptscriptstyle -1}_{\scriptscriptstyle 2} \left(g(\theta,x)\right)^{\scriptscriptstyle \parallel}
\end{equation}
from the orthogonal decomposition of $g = g(\theta,x)$,
\begin{equation} \label{eq:orthogonal}
 g^{\scriptscriptstyle \parallel}\,=\, \mathrm{tr}\left(\theta^{\scriptscriptstyle-1}g\right)\frac{\theta}{m} \hspace{1cm}
g^{\scriptscriptstyle \perp}\,=\, g - g^{\scriptscriptstyle \parallel}
\end{equation}
\end{subequations} 
\indent Figures \ref{fig1} and \ref{fig2} below display numerical results from an application to Kotz-type distributions, which correspond to $h(t) \!=\! -\frac{t^{s}}{\mathstrut 2}$ in (\ref{eq:ecd}) and $\varphi = s^{\scriptscriptstyle 2}\frac{m}{2s}\left(\frac{m}{2s}+1\right)$ in (\ref{eq:varphi})~\cite{kotz}\cite{berkane}. These figures were generated from $10^3$ Monte Carlo runs of the algorithm defined by (\ref{eq:algorithm}) and (\ref{subeq:algoecd}), with random initialisation, for the specific values $s = 4$ and $m = 7$. Essentially the same numerical results could be observed for any $s \leq 9$ and $m \leq 50$. 

Figure \ref{fig1} confirms the fast non-asymptotic  rate of convergence (\ref{eq:ratel2}), stated in (i) of Proposition \ref{prop:information}. On a log-log scale, it shows the empirical mean $\mathbb{E}_{\scriptscriptstyle{\mathrm{MC}}}\,d^{\scriptscriptstyle\, 2}(\theta_{\scriptscriptstyle n},\theta^*)$ over Monte Carlo runs, as a function of $n$. This decreases with a constant negative slope equal to $-1$, starting roughly at $\log n = 4$. Here, the Riemannian distance $d(\theta_{\scriptscriptstyle n},\theta^*)$ induced by the information metric (\ref{subeq:infometric}) is given by~\cite{cyrus}
\begin{equation} \label{eq:infodistance}
 d^{\scriptscriptstyle\,2}(\theta,\theta^*) \,=\, 
I_{\scriptscriptstyle 1}\,\mathrm{tr}\left[\log\left(\theta^{\scriptscriptstyle -1}\theta^*\right)\right]^{\scriptscriptstyle 2}\,+\,I_{\scriptscriptstyle 2}\,\mathrm{tr}^{\scriptscriptstyle 2}\left[\log\left(\theta^{\scriptscriptstyle -1}\theta^*\right)\right]
\hspace{1cm} \theta\,,\theta^* \in \Theta
\end{equation} 
where $\log$ denotes the symmetric matrix logarithm~\cite{higham}. Figure \ref{fig2} confirms the asymptotic efficiency of the recursive estimates $\theta_{\scriptscriptstyle n\,}$, stated in (iii) of Proposition \ref{prop:information}, using Remark \ref{rk:chi2}. It shows a kernel density estimate of $n\,d^{\scriptscriptstyle\,2}(\theta_{\scriptscriptstyle n\,},\theta^*)$ where $n = 10^5$ (solid blue curve). This agrees with a $\chi^{\scriptscriptstyle 2}$-distribution with $28$ degrees of freedom (dotted red curve), where $d = 28$ is indeed the dimension of the parameter space $\mathcal{P}_{\scriptscriptstyle m}$ for $m = 7$. 

\begin{figure}[!b]
  \centering
  \begin{minipage}[b]{0.4\textwidth}
    \includegraphics[width=6cm]{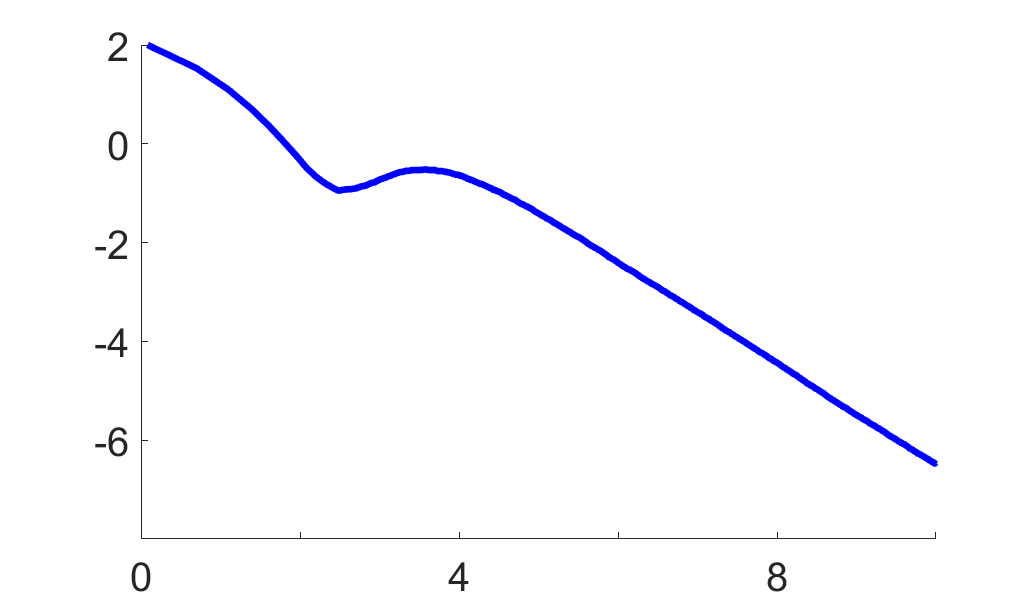}
    \caption{{\small fast non-asymptotic rate of convergence}}
\label{fig1}
  \end{minipage}
  \hfill
  \begin{minipage}[b]{0.4\textwidth}
    \includegraphics[width=6cm]{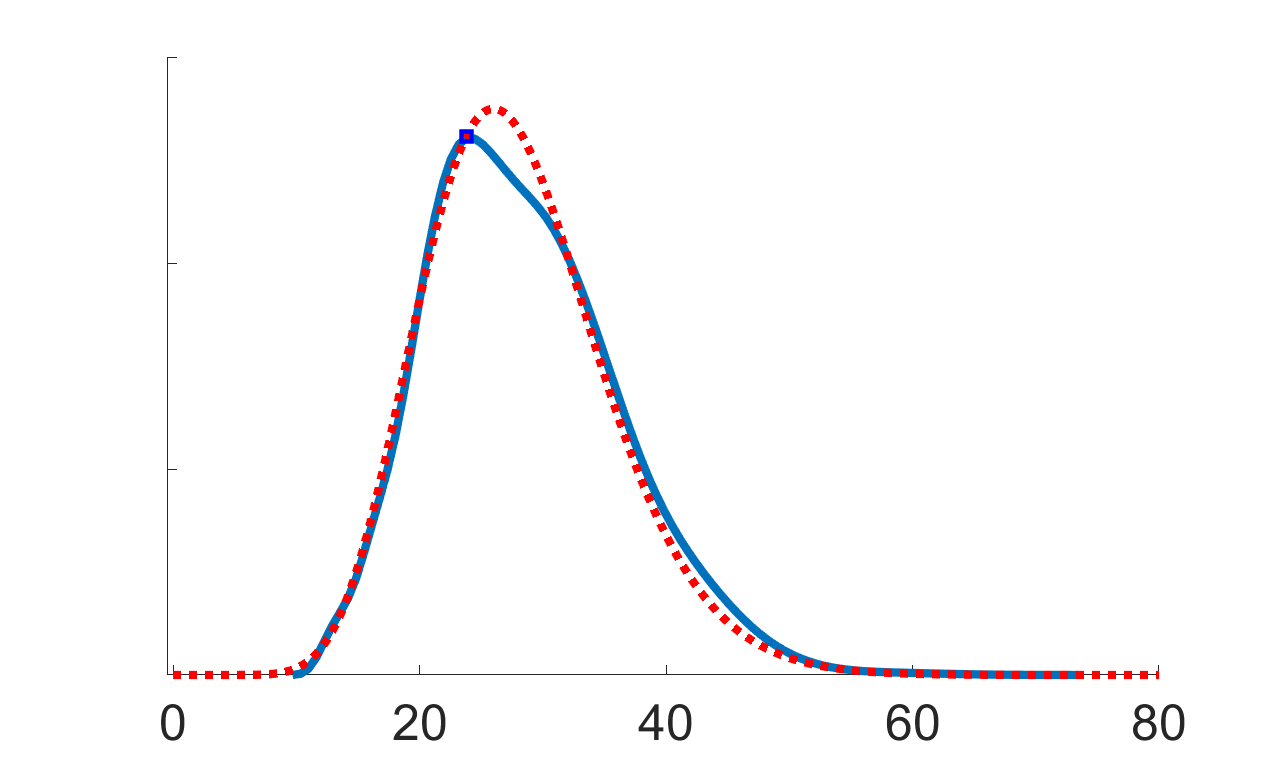}
    \caption{{\small asymptotic efficiency (optimal rate of convergence)}}
\label{fig2}
  \end{minipage}
\end{figure}





\vfill
\pagebreak

\section{Proofs of main results} \label{sec:proofs}
\subsection{Proof of Proposition \ref{prop:as}} \label{subsec:proofas} the proof is a generalisation of the original proof in~\cite{bonnabel}, itself modeled on the proof for the Euclidean case in~\cite{bottou}. Throughout the following, let $\mathcal{X}_{\scriptscriptstyle n}$ be the $\sigma$-field generated by $x_{\scriptscriptstyle 1}\,,\ldots,\,x_{\scriptscriptstyle n}$\!~\cite{shiryayev}. Recall that $(x_{\scriptscriptstyle n}\,;n=1,2,\ldots)$ are i.i.d. with distribution $P_{\scriptscriptstyle \theta^*\,}$. Therefore, by (\ref{eq:algorithm}), $\theta_{\scriptscriptstyle n}$ is $\mathcal{X}_{\scriptscriptstyle n}$-measurable and $x_{\scriptscriptstyle n+1}$ is independent from $\mathcal{X}_{\scriptscriptstyle n\,}$. Thus, using elementary properties of conditional expectation~\cite{shiryayev}, 
\begin{subequations} \label{eq:moments}
\begin{eqnarray}
\label{eq:moments1} \mathbb{E}\left[u(\theta_{\scriptscriptstyle n},x_{\scriptscriptstyle n+1})\middle|\mathcal{X}_{\scriptscriptstyle n}\right] = -D(\theta_{\scriptscriptstyle n}) \\[0.1cm]
\label{eq:moments2} \mathbb{E}\left[\Vert u(\theta_{\scriptscriptstyle n},x_{\scriptscriptstyle n+1})\Vert^2 \middle|\mathcal{X}_{\scriptscriptstyle n}\right] = V(\theta_{\scriptscriptstyle n})
\end{eqnarray}
\end{subequations}
where (\ref{eq:moments1}) follows from (\ref{eq:gradient}), and (\ref{eq:moments2}) from (u1). Let $L$ be a Lipschitz constant for $\nabla D(\theta)$, and $C$ be an upper bound on $V(\theta)$, for $\theta \in \Theta^*$. The following inequality is now proved, for any positive integer $n$,
\begin{equation} \label{eq:rsinequality}
  \mathbb{E}\left[D(\theta_{\scriptscriptstyle n+1}) - D(\theta_{\scriptscriptstyle n})\middle|\mathcal{X}_{\scriptscriptstyle n}\right] \,\leq \gamma^{\scriptscriptstyle 2}_{\scriptscriptstyle n+1}\,LC - \gamma_{\scriptscriptstyle n+1}\Vert \nabla D(\theta_{\scriptscriptstyle n})\Vert^2
\end{equation}
once this is done, Proposition \ref{prop:as} is obtained by applying the Robbins-Siegmund theorem~\cite{duflo}. \\[0.1cm]
\textit{Proof of (\ref{eq:rsinequality})}\,: let $c(t)$ be the geodesic connecting $\theta_{\scriptscriptstyle n}$ to $\theta_{\scriptscriptstyle n+1}$ with equation
\begin{subequations} 
\begin{equation} \label{eq:rsproof1}
   c(t) = \mathrm{Exp}_{\scriptscriptstyle \theta_{\scriptscriptstyle n}}\!\left(t\gamma_{\scriptscriptstyle n+1}u(\theta_{\scriptscriptstyle n},x_{\scriptscriptstyle n+1})\right)
\end{equation}
From the fundamental theorem of calculus,
\begin{equation} \label{eq:rsproof2}
  D(\theta_{\scriptscriptstyle n+1}) - D(\theta_{\scriptscriptstyle n})\,=\, \gamma_{\scriptscriptstyle n+1}\,\langle u(\theta_{\scriptscriptstyle n},x_{\scriptscriptstyle n+1}),\nabla D(\theta_{\scriptscriptstyle n})\rangle \,+\, \gamma_{\scriptscriptstyle n+1}\,\int^1_0\left[\langle \dot{c},\nabla D\rangle_{\scriptscriptstyle c(t)} - \langle \dot{c},\nabla D\rangle_{\scriptscriptstyle c(0)}\right]\,dt
\end{equation}
Since the recursive estimates $\theta_{\scriptscriptstyle n}$ are stable, $\theta_{\scriptscriptstyle n}$ and $\theta_{\scriptscriptstyle n+1}$ both lie in $\Theta^*$. Since $\Theta^*$ is convex, the whole geodesic $c(t)$ lies in $\Theta^*$. Then, since $\nabla D(\theta)$ is Lipschitz on $\Theta^*$, it follows from (\ref{eq:rsproof2}),
\begin{equation} \label{eq:rsproof3}
D(\theta_{\scriptscriptstyle n+1}) - D(\theta_{\scriptscriptstyle n})\,\leq \gamma_{\scriptscriptstyle n+1}\,\langle u(\theta_{\scriptscriptstyle n},x_{\scriptscriptstyle n+1}),\nabla D(\theta_{\scriptscriptstyle n})\rangle \,+\, \gamma^{\scriptscriptstyle 2}_{\scriptscriptstyle n+1}\,L\Vert u(\theta_{\scriptscriptstyle n},x_{\scriptscriptstyle n+1})\Vert^2
\end{equation} 
Taking conditional expectations in this inequality, and using (\ref{eq:moments1}) and (\ref{eq:moments2}),
\begin{equation} \label{eq:rsproof4}
   \mathbb{E}\left[D(\theta_{\scriptscriptstyle n+1}) - D(\theta_{\scriptscriptstyle n})\middle|\mathcal{X}_{\scriptscriptstyle n}\right] \,\leq \,- \gamma_{\scriptscriptstyle n+1}\Vert \nabla D(\theta_{\scriptscriptstyle n})\Vert^2 \,+\,
\gamma^{\scriptscriptstyle 2}_{\scriptscriptstyle n+1}\,LV(\theta_{\scriptscriptstyle n})
\end{equation}
\end{subequations}
so (\ref{eq:rsinequality}) follows since (u1) guarantees $V(\theta_{\scriptscriptstyle n}) \leq C$. \hfill$\blacksquare$ \\[0.1cm]
\textit{Conclusion}\,: by the Robbins-Siegmund theorem, inequality (\ref{eq:rsinequality}) implies that, almost surely,
\begin{subequations}
\begin{equation} \label{eq:proofas1}
  \lim D(\theta_{\scriptscriptstyle n}) =D_{\scriptscriptstyle \infty} \,<\infty \hspace{0.25cm}and\hspace{0.25cm} \sum^{\scriptscriptstyle\infty}_{\scriptscriptstyle n=1}\,\gamma_{\scriptscriptstyle n+1}\,\Vert \nabla D(\theta_{\scriptscriptstyle n})\Vert^2 \, < \infty
\end{equation} 
In particular, from the first condition in (\ref{eq:stepsize}), convergence of the sum in (\ref{eq:proofas1}) implies
\begin{equation} \label{eq:proofas2}
  \lim\,\Vert\nabla D(\theta_{\scriptscriptstyle n})\Vert \,= 0 \hspace{1cm} \text{almost surely}
\end{equation}
\end{subequations}
Now, since the sequence of recursive estimates $\theta_{\scriptscriptstyle n}$ lies in the compact set $\Theta^*$, it has at least one point of accumulation in this set, say $\theta_{*\,}$. If $\theta_{\scriptscriptstyle n(k)}$ is a subsequence of $\theta_{\scriptscriptstyle n\,}$, converging to $\theta_{*\,}$,
$$
\Vert \nabla D(\theta_*)\Vert \,=\, \lim\,\Vert\nabla D(\theta_{\scriptscriptstyle n(k)})\Vert\,=\,   \lim\,\Vert\nabla D(\theta_{\scriptscriptstyle n})\Vert \,= 0 \hspace{1cm} \text{almost surely}
$$
where the third equality follows from (\ref{eq:proofas2}). This means that $\theta_*$ is a stationary point of $D(\theta)$ in $\Theta^*$. Thus, (d1) implies $\theta_* = \theta^*$ is the unique point of accumulation of $\theta_{\scriptscriptstyle n\,}$. In other words, $\lim \theta_{\scriptscriptstyle n} = \theta^*$ almost surely. \hfill $\blacksquare$ 
\subsection{Proof of Proposition \ref{prop:ratel2}} the proof is modeled on the proofs for the Euclidean case, given in~\cite{nev}\cite{benv}. It relies on the following geometric Lemmas \ref{lemma:grad} and \ref{lemma:trigo}. Lemma \ref{lemma:grad} will be proved in Appendix \ref{sec:geometric}. On the other hand, Lemma \ref{lemma:trigo} is the same as the trigonometric distance bound of~\cite{sra}. For Lemma \ref{lemma:grad}, recall that $\lambda > 0$ denotes the smallest eigenvalue of the matrix $H$ defined in (\ref{eq:hess}).
\begin{subequations}
\begin{lemma} \label{lemma:grad}
  for any $\mu < \lambda$, there exists a neighborhood $\bar{\Theta}^*$ of $\theta^*$, contained in $\Theta^*$, with
\begin{equation} \label{eq:lemmgrad}
   \langle\mathrm{Exp}^{\scriptscriptstyle -1}_{\scriptscriptstyle \theta}(\theta^*),\nabla D(\theta)\rangle \,\leq\, -\mu\,d^{\scriptscriptstyle\,2}(\theta,\theta^*) \hspace{1cm}\text{for } \theta \in \bar{\Theta}^*
\end{equation}
\end{lemma}
\begin{lemma} \label{lemma:trigo}
let $-\kappa^{\scriptscriptstyle 2}$ be a lower bound on the sectional curvature of $\Theta$ in $\Theta^*$, and $C_{\scriptscriptstyle\kappa} = R\kappa\coth(R\kappa)$ where $R$ is the diameter of $\Theta^*$. For $\tau,\theta \in \Theta^*$, where $\tau = \mathrm{Exp}_{\scriptscriptstyle \theta}(u)$,
\begin{equation} \label{eq:trigo}
  d^{\scriptscriptstyle\,2}(\tau,\theta^*) \,\leq\, d^{\scriptscriptstyle\,2}(\theta,\theta^*) - 2\,\langle \mathrm{Exp}^{\scriptscriptstyle -1}_{\scriptscriptstyle \theta}(\theta^*),u\rangle + C_{\scriptscriptstyle\kappa}\Vert u\Vert^2
\end{equation}
\end{lemma}
\end{subequations}
\noindent \textit{Proof of (\ref{eq:ratel2})}\,: let $\gamma_{\scriptscriptstyle n} = \frac{a}{n}$ with $2\lambda a > 2\mu a > 1$ for some $\mu < \lambda$, and let $\bar{\Theta}^*$ be the neighborhood corresponding to $\mu$ in Lemma \ref{lemma:grad}. By Proposition \ref{prop:as}, the $\theta_{\scriptscriptstyle n}$ converge to $\theta^*$ almost surely. Without loss of generality, it can be assumed that all the $\theta_n$ lie in $\bar{\Theta}^*$, almost surely. Then, (\ref{eq:algorithm}) and Lemma \ref{lemma:trigo} imply, for any positive integer $n$,
\begin{subequations}
\begin{equation} \label{eq:trigalgo1}
  d^{\scriptscriptstyle\,2}(\theta_{\scriptscriptstyle n+1},\theta^*) \,\leq\, d^{\scriptscriptstyle\,2}(\theta_{\scriptscriptstyle n},\theta^*) - 2\gamma_{\scriptscriptstyle n+1}\,\langle \mathrm{Exp}^{\scriptscriptstyle -1}_{\scriptscriptstyle \theta_{\scriptscriptstyle n}}(\theta^*),u(\theta_{\scriptscriptstyle n},x_{\scriptscriptstyle n+1})\rangle + \gamma^{\scriptscriptstyle 2}_{\scriptscriptstyle n+1}\,C_{\scriptscriptstyle\kappa}\Vert u(\theta_{\scriptscriptstyle n},x_{\scriptscriptstyle n+1})\Vert^2
\end{equation} 
Indeed, this follows by replacing $\tau = \theta_{\scriptscriptstyle n+1}$ and $\theta = \theta_{\scriptscriptstyle n}$ in (\ref{eq:trigo}). Taking conditional expectations in (\ref{eq:trigalgo1}), and using (\ref{eq:moments1}) and (\ref{eq:moments2}),
$$
 \mathbb{E}\left[d^{\scriptscriptstyle\,2}(\theta_{\scriptscriptstyle n+1},\theta^*)\middle|\mathcal{X}_{\scriptscriptstyle n}\right]\,\leq\, 
d^{\scriptscriptstyle\,2}(\theta_{\scriptscriptstyle n},\theta^*) + 2\gamma_{\scriptscriptstyle n+1}\,\langle \mathrm{Exp}^{\scriptscriptstyle -1}_{\scriptscriptstyle \theta_{\scriptscriptstyle n}}(\theta^*),\nabla D(\theta_{\scriptscriptstyle n})\rangle + \gamma^{\scriptscriptstyle 2}_{\scriptscriptstyle n+1}\,C_{\scriptscriptstyle\kappa}V(\theta_{\scriptscriptstyle n})
$$
Then, by (u1) and (\ref{eq:lemmgrad}) of Lemma \ref{lemma:grad},
\begin{equation} \label{eq:trigalgo3}
\mathbb{E}\left[d^{\scriptscriptstyle\,2}(\theta_{\scriptscriptstyle n+1},\theta^*)\middle|\mathcal{X}_{\scriptscriptstyle n}\right]\,\leq\,
d^{\scriptscriptstyle\,2}(\theta_{\scriptscriptstyle n},\theta^*)(1-2\gamma_{\scriptscriptstyle n+1}\mu) + 
\gamma^{\scriptscriptstyle 2}_{\scriptscriptstyle n+1}\,C_{\scriptscriptstyle\kappa}C
\end{equation}
where $C$ is an upper bound on $V(\theta)$, for $\theta \in \Theta^*$. By further taking expectations 
\begin{equation} \label{eq:trigalgo4}
   \mathbb{E}\,d^{\scriptscriptstyle\,2}(\theta_{\scriptscriptstyle n+1},\theta^*)\,\leq\,
\mathbb{E}\,d^{\scriptscriptstyle\,2}(\theta_{\scriptscriptstyle n},\theta^*)(1-2\gamma_{\scriptscriptstyle n+1}\mu) + 
\gamma^{\scriptscriptstyle 2}_{\scriptscriptstyle n+1}\,C_{\scriptscriptstyle\kappa}C
\end{equation}
\end{subequations}
Using (\ref{eq:trigalgo4}), the proof reduces to an elementary reasoning by recurrence. Indeed, replacing $\gamma_{\scriptscriptstyle n} = \frac{a}{n}$ into (\ref{eq:trigalgo4}), it follows that 
\begin{subequations}
\begin{equation} \label{eq:recurrence1}
   \mathbb{E}\,d^{\scriptscriptstyle\,2}(\theta_{\scriptscriptstyle n+1},\theta^*)\,\leq\,
\mathbb{E}\,d^{\scriptscriptstyle\,2}(\theta_{\scriptscriptstyle n},\theta^*)\left(1-\frac{2\mu a}{ n+1}\right) + 
\frac{a^{\scriptscriptstyle 2}C_{\scriptscriptstyle\kappa}C}{(n+1)^{\scriptscriptstyle 2}}   
\end{equation}
On the other hand, if $b(n) = \frac{b}{n}$ where $b > a^{\scriptscriptstyle 2}C_{\scriptscriptstyle\kappa}C\,(2\mu a -1)^{\scriptscriptstyle -1}$, then
\begin{equation} \label{eq:recurrence2}
  b(n+1) \geq b(n) \left(1-\frac{2\mu a}{ n+1}\right) + 
\frac{a^{\scriptscriptstyle 2}C_{\scriptscriptstyle\kappa}C}{(n+1)^{\scriptscriptstyle 2}} \\[0.1cm]
\end{equation}
\end{subequations}
Let $b$ be sufficiently large, so (\ref{eq:recurrence2}) is verified and $\mathbb{E}\,d^{\scriptscriptstyle\,2}(\theta_{\scriptscriptstyle n_{\scriptscriptstyle o}},\theta^*) \leq b(n_{\scriptscriptstyle o})$ for some $n_{\scriptscriptstyle o\,}$. Then, by recurrence, using (\ref{eq:recurrence1}) and (\ref{eq:recurrence2}), one also has that $\mathbb{E}\,d^{\scriptscriptstyle\,2}(\theta_{\scriptscriptstyle n\,},\theta^*) \leq b(n)$ for all $n \geq n_{\scriptscriptstyle o\,}$. In other words, (\ref{eq:ratel2}) holds true. \hfill $\blacksquare$ \vfill
\pagebreak
\subsection{Proof of Proposition \ref{prop:rateas}} 
the proof is modeled on the proof for the Euclidean case in~\cite{nev}. To begin, let $W_{\scriptscriptstyle n}$ be the stochastic process given by
\begin{subequations}
\begin{equation} \label{eq:proofrateas1}
  W_{\scriptscriptstyle n} \,=\, n^{\scriptscriptstyle p}\,d^{\scriptscriptstyle\,2}(\theta_{\scriptscriptstyle n},\theta^*) + n^{\scriptscriptstyle -q} \hspace{1cm} \text{ where } q \in (0,1-p)
\end{equation}
The idea is to show that this process is a positive supermartingale, for sufficiently large $n$. By the supermartingale convergence theorem~\cite{shiryayev}, it then follows that $W_{\scriptscriptstyle n}$ converges to a finite limit, almost surely. In particular, this implies 
 \begin{equation} \label{eq:proofrateas2}
  \lim n^{\scriptscriptstyle p}\,d^{\scriptscriptstyle\,2}(\theta_{\scriptscriptstyle n},\theta^*) \,=\, \ell_{\scriptscriptstyle p} < \infty \hspace{1cm} \text{almost surely}
 \end{equation}
Then, $\ell_{\scriptscriptstyle p}$ must be equal to zero, since $p$ is arbitrary in the interval $(0,1)$. Precisely, for any $\varepsilon \in (0,1-p)$,
$$
\ell_{\scriptscriptstyle p} \,=\, \lim n^{\scriptscriptstyle p}\,d^{\scriptscriptstyle\,2}(\theta_{\scriptscriptstyle n},\theta^*)
\,=\, \lim n^{\scriptscriptstyle -\varepsilon }n^{\scriptscriptstyle p+\varepsilon}\,d^{\scriptscriptstyle\,2}(\theta_{\scriptscriptstyle n},\theta^*)
\,=\,\left(\lim n^{\scriptscriptstyle -\varepsilon }\right)\,\ell_{\scriptscriptstyle p+\varepsilon} \,=\, 0
$$
\end{subequations}
It remains to show that $W_{\scriptscriptstyle n}$ is a supermartingale, for sufficiently large $n$. To do so, note that by (\ref{eq:trigalgo3}) from the proof of Proposition \ref{prop:ratel2},
$$
\mathbb{E}\left[W_{\scriptscriptstyle n+1}-W_{\scriptscriptstyle n}\middle|\mathcal{X}_{\scriptscriptstyle n}\right] \leq\, d^{\scriptscriptstyle\,2}(\theta_{\scriptscriptstyle n},\theta^*)\,\frac{p-2\mu a}{(n+1)^{\scriptscriptstyle1-p}} \,+\, \frac{a^{\scriptscriptstyle 2}C_{\scriptscriptstyle\kappa}C}{(n+1)^{\scriptscriptstyle 2-p}} \,-\, \frac{q}{(n+1)^{\scriptscriptstyle q+1}}   
$$
Here, the first term on the right-hand side is negative, since $2\mu a > 1 > p$. Moreover, the third term dominates the second one for sufficiently large $n$, since $ q < 1- p$. Thus, for sufficiently large $n$, the right-hand side is negative, and $W_{\scriptscriptstyle n}$ is a supermartingale.\hfill $\blacksquare$
\subsection{Proof of Proposition \ref{prop:normality}} the proof relies on the following geometric Lemmas \ref{lemma:linearalgo} and \ref{lemma:linearfield}, which are used to linearise Algorithm (\ref{eq:algorithm}), in terms of the normal coordinates $\theta^{\scriptscriptstyle\,\alpha}$. This idea of linearisation in terms of local coordinates also plays a central role in~\cite{flamm}. 
\begin{subequations}
\begin{lemma} \label{lemma:linearalgo}
let $\theta_{\scriptscriptstyle n\,},\theta_{\scriptscriptstyle n+1}$ be given by (\ref{eq:algorithm}) with $\gamma_{\scriptscriptstyle n} = \frac{a}{n\,}$. Then, in a system of normal coordinates with origin at $\theta^*$,
\begin{equation} \label{eq:linearalgo}
  \theta^{\scriptscriptstyle\,\alpha}_{\scriptscriptstyle n+1}\,=\,   \theta^{\scriptscriptstyle\,\alpha}_{\scriptscriptstyle n}+\gamma^{\phantom{\scriptscriptstyle 2}}_{\scriptscriptstyle n+1}\,u^{\scriptscriptstyle \alpha}_{\scriptscriptstyle n+1} + \gamma^{\scriptscriptstyle 2}_{\scriptscriptstyle n+1}\,\pi^{\scriptscriptstyle \alpha}_{\scriptscriptstyle n+1} \hspace{1cm} \mathbb{E}\left|\pi^{\scriptscriptstyle \alpha}_{\scriptscriptstyle n+1}\right| = O(n^{\scriptscriptstyle -1/2})
\end{equation}
where $u^{\scriptscriptstyle \alpha}_{\scriptscriptstyle n+1}$ are the components of $u(\theta_{\scriptscriptstyle n},x_{\scriptscriptstyle n+1})$.
\end{lemma}
\begin{lemma} \label{lemma:linearfield}
  let $v_{\scriptscriptstyle n} = \nabla D(\theta_{\scriptscriptstyle n})\,$. Then, in a system of normal coordinates with origin at $\theta^*$,
\begin{equation} \label{eq:linearfield}
  v^{\scriptscriptstyle\, \alpha}_{\scriptscriptstyle n}\,=\, H^{\phantom{\scriptscriptstyle 2}}_{\scriptscriptstyle\alpha\beta}\,\theta^{\scriptscriptstyle\,\beta}_{\scriptscriptstyle n}\,+\, \rho^{\scriptscriptstyle\alpha}_{\scriptscriptstyle n} \hspace{1cm} \rho^{\scriptscriptstyle\alpha}_{\scriptscriptstyle n} = o\left( d(\theta_{\scriptscriptstyle n},\theta^*)\right)
\end{equation}
where $v^{\scriptscriptstyle\, \alpha}_{\scriptscriptstyle n}$ are the components of $v_{\scriptscriptstyle n}$ and the $H_{\scriptscriptstyle\alpha\beta}$ were defined in (\ref{eq:hess}). 
\end{lemma}
\end{subequations}
\noindent \textit{Linearisation of (\ref{eq:algorithm})}\,: let $u(\theta_{\scriptscriptstyle n},x_{\scriptscriptstyle n+1}) = -v_{\scriptscriptstyle n}+w_{\scriptscriptstyle n+1\,}$. Then, it follows from (\ref{eq:linearalgo}) and (\ref{eq:linearfield}),
\begin{subequations} 
\begin{equation} \label{eq:linearisation1}
   \theta^{\scriptscriptstyle\,\alpha}_{\scriptscriptstyle n+1} \,=\,    \theta^{\scriptscriptstyle\,\alpha}_{\scriptscriptstyle n} \,-\,
\gamma^{\phantom{\scriptscriptstyle 2}}_{\scriptscriptstyle n+1}\,H^{\phantom{\scriptscriptstyle 2}}_{\scriptscriptstyle\alpha\beta}\,\theta^{\scriptscriptstyle\,\beta}_{\scriptscriptstyle n}\,-\,\gamma^{\phantom{\scriptscriptstyle 2}}_{\scriptscriptstyle n+1}\,\rho^{\scriptscriptstyle\alpha}_{\scriptscriptstyle n}\,+\,\gamma^{\phantom{\scriptscriptstyle 2}}_{\scriptscriptstyle n+1}\,w^{\scriptscriptstyle\,\alpha}_{\scriptscriptstyle n+1}\,+\,\gamma^{\scriptscriptstyle 2}_{\scriptscriptstyle n+1}\,\pi^{\scriptscriptstyle \alpha}_{\scriptscriptstyle n+1}
\end{equation}
Denote the re-scaled coordinates $n^{\scriptscriptstyle 1/2}\theta^{\scriptscriptstyle\,\alpha}_{\scriptscriptstyle n}$ by $\eta^{\scriptscriptstyle\alpha}_{\scriptscriptstyle n\,}$, and recall $\gamma_{\scriptscriptstyle n} = \frac{a}{n\,}$. Then, using the estimate $(n+1)^{\scriptscriptstyle 1/2} = n^{\scriptscriptstyle 1/2}(1+(2n)^{\scriptscriptstyle -1}+O(n^{\scriptscriptstyle -2}))$, it follows from (\ref{eq:linearisation1}) that 
\begin{equation} \label{eq:linearisation2}
  \eta^{\scriptscriptstyle\alpha}_{\scriptscriptstyle n+1}\,=\,   \eta^{\scriptscriptstyle\alpha}_{\scriptscriptstyle n} + \frac{A_{\scriptscriptstyle\alpha\beta}}{n+1}\,\eta^{\scriptscriptstyle\beta}_{\scriptscriptstyle n}\,+\, \frac{a}{(n+1)^{\scriptscriptstyle 1/2}}\,\left[ B^{\phantom{\scriptscriptstyle2}}_{\scriptscriptstyle\alpha\beta}\,\theta^{\scriptscriptstyle\,\beta}_{\scriptscriptstyle n} \,- \rho^{\scriptscriptstyle\alpha}_{\scriptscriptstyle n}\,+\,w^{\scriptscriptstyle\,\alpha}_{\scriptscriptstyle n+1}\,+\,\frac{a\pi^{\scriptscriptstyle\alpha}_{\scriptscriptstyle n+1}}{n+1}\,\right] \\[0.1cm]
\end{equation} 
where $A_{\scriptscriptstyle \alpha\beta} = \frac{1}{2}\delta_{\alpha\beta}-aH_{\scriptscriptstyle\alpha\beta}$ and $B_{\scriptscriptstyle\alpha\beta} = O(n^{\scriptscriptstyle-1})$. Equation (\ref{eq:linearisation2}) is a first-order, inhomogeneous, linear difference equation, for the ``vector" $\eta^{\phantom{\scriptscriptstyle\,2}}_{\scriptscriptstyle n}$ of components $\eta^{\scriptscriptstyle\alpha}_{\scriptscriptstyle n\,}$. \hfill$\blacksquare$ 
\end{subequations}
\vfill
\pagebreak
\noindent \textit{Study of equation (\ref{eq:linearisation2})}\,: switching to vector-matrix notation, equation (\ref{eq:linearisation2}) is of the general form
\begin{subequations}
  \begin{equation} \label{eq:linearisation3}
     \eta^{\phantom{\scriptscriptstyle\,2}}_{\scriptscriptstyle n+1} \,=\, \left(I \,+\,\frac{A}{n+1}\right)\,\eta^{\phantom{\scriptscriptstyle\,2}}_{\scriptscriptstyle n} \,+\, \frac{a\,\xi^{\phantom{\scriptscriptstyle\,2}}_{\scriptscriptstyle n+1}}{(n+1)^{\scriptscriptstyle 1/2}}
  \end{equation}
where $I$ denotes the identity matrix, $A$ has matrix elements $A_{\scriptscriptstyle \alpha\beta\,}$, and $\left(\xi_{\scriptscriptstyle n}\right)$ is a sequence of inputs. The general solution of this equation is~\cite{nev}\cite{kailath}
\begin{equation} \label{eq:transition1}
   \eta^{\phantom{\scriptscriptstyle\,2}}_{\scriptscriptstyle n}\,=\, A^{\phantom{\scriptscriptstyle\,2}}_{\scriptscriptstyle n,m}\,\eta^{\phantom{\scriptscriptstyle\,2}}_{\scriptscriptstyle m}\,+\, \sum^{\scriptscriptstyle n}_{\scriptscriptstyle k = m+1}\,A^{\phantom{\scriptscriptstyle\,2}}_{\scriptscriptstyle n,k}\,\frac{a\,\xi^{\phantom{\scriptscriptstyle\,2}}_{\scriptscriptstyle k}}{k^{\scriptscriptstyle 1/2}} \hspace{1cm} \text{for }\,\, n \geq m
\end{equation} 
where the transition matrix $A_{\scriptscriptstyle n,k}$ is given by
\begin{equation} \label{eq:transitions2}
   A^{\phantom{\scriptscriptstyle\,2}}_{\scriptscriptstyle n,k} \,=\,\prod^{\scriptscriptstyle n}_{\scriptscriptstyle j = k+1}\,\left(I+\frac{A}{j}\right) \hspace{1cm}
A^{\phantom{\scriptscriptstyle\,2}}_{\scriptscriptstyle n,n} = I
\end{equation}
Since $2\lambda a > 1$, the matrix $A$ is stable. This can be used to show that~\cite{nev}\cite{kailath} 
\begin{equation} \label{eq:nev}
  q > \frac{1}{2} \,\text{ and }\,\mathbb{E}\left|\xi^{\phantom{\scriptscriptstyle\,2}}_{\scriptscriptstyle n}\right| = O(n^{\scriptscriptstyle -q}) \,\,\Longrightarrow\,\, \lim\eta^{\phantom{\scriptscriptstyle\,2}}_{\scriptscriptstyle n} \,=\, 0\,\text{ in probability}
\end{equation}
where $|\xi_{\scriptscriptstyle n}|$ denotes the Euclidean vector norm. Then, it follows from (\ref{eq:nev}) that  $\eta^{\phantom{\scriptscriptstyle\,2}}_{\scriptscriptstyle n}$ converges to zero in probability, in each one of the three cases
\end{subequations}
$$
\xi^{\scriptscriptstyle\alpha}_{\scriptscriptstyle n+1} \,=\, B^{\phantom{\scriptscriptstyle\,2}}_{\scriptscriptstyle \alpha\beta}\,\theta^{\scriptscriptstyle\,\beta}_{\scriptscriptstyle n}\hspace{0.25cm};\hspace{0.25cm}
\xi^{\scriptscriptstyle\alpha}_{\scriptscriptstyle n+1} \,=\, \rho^{\scriptscriptstyle\alpha}_{\scriptscriptstyle n}\hspace{0.25cm};\hspace{0.25cm}
\xi^{\scriptscriptstyle\alpha}_{\scriptscriptstyle n+1} \,=\, \frac{\pi^{\scriptscriptstyle\alpha}_{\scriptscriptstyle n+1}}{n+1}
$$
Indeed, in the first two cases, the condition required in (\ref{eq:nev}) can be verified using (\ref{eq:ratel2}), whereas in the third case, it follows immediately from the estimate of $\mathbb{E}|\pi^{\scriptscriptstyle\alpha}_{\scriptscriptstyle n+1}|$ in (\ref{eq:linearalgo}).
\hfill$\blacksquare$ \\[0.1cm]
\noindent \textit{Conclusion}\,: by linearity of (\ref{eq:linearisation2}),  it is enough to consider the case $\xi^{\scriptscriptstyle\alpha}_{\scriptscriptstyle n+1} = w^{\scriptscriptstyle\,\alpha}_{\scriptscriptstyle n+1}$ in (\ref{eq:linearisation3}). Then, according to (\ref{eq:transition1}), $\eta^{\phantom{\scriptscriptstyle\,2}}_{\scriptscriptstyle n}$ has the same limit distribution as the sums
\begin{equation} \label{eq:sum}
  \tilde{\eta}^{\phantom{\scriptscriptstyle\,2}}_{\scriptscriptstyle n} \,=\, \sum^{\scriptscriptstyle n}_{\scriptscriptstyle k=1}\,
A^{\phantom{\scriptscriptstyle\,2}}_{\scriptscriptstyle n,k}\,\frac{aw_{\scriptscriptstyle k}}{k^{\scriptscriptstyle 1/2}}
\end{equation}
\begin{subequations} \label{subeq:clt}
By (\ref{eq:moments}), $(w_{\scriptscriptstyle k})$ is a sequence of square-integrable martingale differences. Therefore, to conclude that the limit distribution of $\tilde{\eta}^{\phantom{\scriptscriptstyle\,2}}_{\scriptscriptstyle n}$ is a centred $d$-variate normal distribution, with covariance matrix $\Sigma$ given by (\ref{eq:lyapunov}), it is enough to verify the conditions of the martingale central limit theorem~\cite{martingale},
\begin{equation} \label{eq:clt1}
\lim\max_{\scriptscriptstyle k\leq n}\,\left| A^{\phantom{\scriptscriptstyle\,2}}_{\scriptscriptstyle n,k}\,\frac{aw_{\scriptscriptstyle k}}{k^{\scriptscriptstyle 1/2}}\right| \,=\, 0 \,\text{ in probability}\hspace{0.51cm}
\end{equation}
\begin{equation} \label{eq:clt2}
\sup\,\mathbb{E}\left|\tilde{\eta}^{\phantom{\scriptscriptstyle\,2}}_{\scriptscriptstyle n}\right|^2 \,<\,\infty \hspace{4.2cm}
\end{equation}
\begin{equation} \label{eq:clt3}
\lim \sum^{\scriptscriptstyle n}_{\scriptscriptstyle k=1}\frac{a^{\scriptscriptstyle 2}}{k}\,A^{\phantom{\scriptscriptstyle\,2}}_{\scriptscriptstyle n,k\,}\Sigma^{\phantom{\scriptscriptstyle\,2}}_{\scriptscriptstyle k\,}A^{\phantom{\scriptscriptstyle\,2}}_{\scriptscriptstyle n,k} \,=\,\Sigma \,\text{ in probability}
\end{equation}
\end{subequations}
where $\Sigma^{\phantom{\scriptscriptstyle\,2}}_{\scriptscriptstyle k}$ is the conditional covariance matrix
\begin{equation} \label{eq:sigmak}
  \Sigma^{\phantom{\scriptscriptstyle\,2}}_{\scriptscriptstyle k} \,=\,\mathbb{E}\left[w^{\phantom{\scriptscriptstyle\dagger}}_{\scriptscriptstyle k}w^{\scriptscriptstyle{\dagger}}_{\scriptscriptstyle k\,}\middle|\mathcal{X}_{\scriptscriptstyle k-1}\right]
\end{equation}
Conditions (\ref{subeq:clt}) are verified in Appendix \ref{sec:clt},  which completes the proof. \hfill$\blacksquare$
\subsection{Proof of Proposition \ref{prop:information}} denote $\partial_{\scriptscriptstyle\alpha}=\frac{\partial}{\mathstrut\partial\theta^{\scriptscriptstyle\,\alpha}}$ the coordinate vector fields of the normal coordinates $\theta^{\scriptscriptstyle\,\alpha\,}$. Since $\langle\cdot,\cdot\rangle$ coincides with the information metric of the model $P$, it follows from (\ref{eq:hess}) and (\ref{eq:raofish2}), 
\begin{subequations}
\begin{equation} \label{eq:hid1}
 H_{\scriptscriptstyle\alpha\beta} \,=\, \langle \partial_{\scriptscriptstyle \alpha\,},\partial_{\scriptscriptstyle \beta}\rangle_{\scriptscriptstyle \theta^*}
\end{equation}
However, by the definition of normal coordinates~\cite{petersen}, the $\partial_{\scriptscriptstyle \alpha}$ are orthonormal at $\theta^*$. Therefore,
\begin{equation} \label{eq:hid2}
 H_{\scriptscriptstyle\alpha\beta} \,=\, \delta_{\scriptscriptstyle\alpha\beta}
\end{equation}
\end{subequations}
Thus, the matrix $H$ is equal to the identity matrix, and its smallest eigenvalue is $\lambda = 1$. \\[0.1cm]
\textit{Proof of \emph{(i)}}\,: this follows directly from Propositions \ref{prop:ratel2} and \ref{prop:rateas}. Indeed, since $\lambda = 1$, the conditions of these propositions are verified, as soon as $2a > 1$. Therefore, (\ref{eq:ratel2}) and (\ref{eq:rateas}) hold true. \hfill$\blacksquare$ \\[0.1cm]
\textit{Proof of \emph{(ii)}}\,: this follows from Proposition \ref{prop:normality}. The conditions of this proposition are verified, as soon as $2a > 1$. Therefore, the distribution of the re-scaled coordinates $(n^{\scriptscriptstyle 1/2}\theta^{\scriptscriptstyle\,\alpha}_{\scriptscriptstyle n})$ converges to a centred $d$-variate normal distribution, with covariance matrix $\Sigma$ given by Lyapunov's equation (\ref{eq:lyapunov}). If $a = 1$, then (\ref{eq:hid2}) implies $A_{\scriptscriptstyle\alpha\beta} = - \frac{1}{2}\delta_{\scriptscriptstyle \alpha\beta\,}$, so that Lyapunov's equation (\ref{eq:lyapunov}) reads $\Sigma = \Sigma^*$, as required. \hfill$\blacksquare$ \\[0.1cm]
\indent For the following proof of (iii), the reader may wish to recall that summation convention is used throughout the present work. That is~\cite{petersen}, summation is implicitly  understood over any repeated subscript or superscript from the Greek alphabet, taking the values $1\,,\ldots,\,d\,$.\\[0.1cm]
\textit{Proof of \emph{(iii)}}\,: let $\ell(\theta) = \log L(\theta)$ and assume $u(\theta,x)$ is given by (\ref{eq:score}). Then, by the definition of normal coordinates~\cite{petersen}, the following expression holds
\begin{subequations}
\begin{equation} \label{eq:score1}
   u^{\scriptscriptstyle \alpha}(\theta^*) \,=\, \left.\frac{\partial \ell}{\partial \theta^{\scriptscriptstyle\,\alpha}}\right|_{\scriptscriptstyle \theta^{\scriptscriptstyle \alpha} = 0}
\end{equation}
Replacing this into (\ref{eq:cov}) gives
\begin{equation} \label{eq:score2}
  \Sigma^*_{\scriptscriptstyle\alpha\beta} \,=\, E_{\scriptscriptstyle \theta^*}\! \left[\frac{\partial \ell}{\partial \theta^{\scriptscriptstyle\,\alpha}}\frac{\partial \ell}{\partial \theta^{\scriptscriptstyle\,\beta}}\right]_{\scriptscriptstyle \theta^{\scriptscriptstyle \alpha} = 0} \,=\, 
 - \,E_{\scriptscriptstyle \theta^*}\left.\frac{\partial^{\scriptscriptstyle\, 2}\!\,\ell}{\mathstrut\partial\theta^{\scriptscriptstyle \alpha}\partial\theta^{\scriptscriptstyle \beta}}\right|_{\scriptscriptstyle \theta^{\scriptscriptstyle \alpha} = 0} \,=\, \left.\frac{\partial^{\scriptscriptstyle\, 2}\!\,D}{\mathstrut\partial\theta^{\scriptscriptstyle \alpha}\partial\theta^{\scriptscriptstyle \beta}}\right|_{\scriptscriptstyle \theta^{\scriptscriptstyle \alpha} = 0} 
\end{equation}
\end{subequations}
where the second equality is the so-called Fisher's identity (see~\cite{amari}, Page 28), and the third equality follows from (\ref{eq:kl}) by differentiating under the expectation. Now, by (\ref{eq:hess}) and (\ref{eq:hid2}), $\Sigma^*$ is the identity matrix. 

To show that the recursive estimates $\theta_{\scriptscriptstyle n}$ are asymptotically efficient, let $(\tau^{\scriptscriptstyle\alpha}\,;\alpha = 1,\ldots, d\,)$ be any local coordinates with origin at $\theta^*$ and let $\tau^{\scriptscriptstyle \alpha}_{\scriptscriptstyle n} = \tau^{\scriptscriptstyle \alpha}(\theta_{\scriptscriptstyle n})\,$. From the second-order Taylor expansion of each coordinate function $\tau^{\scriptscriptstyle \alpha}$, it is straightforward to show that
\begin{subequations} 
\begin{equation} \label{eq:proofeff1}  
n^{\scriptscriptstyle 1/2}\tau^{\scriptscriptstyle\alpha}_{\scriptscriptstyle n} \,=\,\left(\frac{\partial \tau^{\scriptscriptstyle\alpha}}{\partial\theta^{\scriptscriptstyle\,\gamma}}\right)_{\!\scriptscriptstyle \theta^*} \!\left(n^{\scriptscriptstyle 1/2}\theta^{\scriptscriptstyle\,\gamma}_{\scriptscriptstyle n\,}\right)\,+\, \sigma^{\scriptscriptstyle\alpha}(\theta_{\scriptscriptstyle n})\left(n^{\scriptscriptstyle 1/2}d^{\scriptscriptstyle\, 2}(\theta_{\scriptscriptstyle n\,},\theta^*)\right)
\end{equation}
where the subscript $\theta^*$ indicates the derivative is evaluated at $\theta^*$, and where $\sigma^{\scriptscriptstyle\alpha}$ is a continuous function in the neighborhood of $\theta^*$. By (\ref{eq:rateas}), the second term in (\ref{eq:proofeff1}) converges to zero almost surely. Therefore, the limit distribution of the re-scaled coordinates $(n^{\scriptscriptstyle 1/2}\tau^{\scriptscriptstyle\alpha}_{\scriptscriptstyle n})$ is the same as that of the first term in (\ref{eq:proofeff1}). By (ii), this is a centred $d$-variate normal distribution with covariance matrix  $\Sigma^{\scriptscriptstyle\tau}$ given by
\begin{equation} \label{eq:proofeff2}
\Sigma^{\scriptscriptstyle\tau}_{\scriptscriptstyle\alpha\beta} \,=\,
\left(\frac{\partial \tau^{\scriptscriptstyle\alpha}}{\partial\theta^{\scriptscriptstyle\,\gamma}}\right)_{\!\scriptscriptstyle \theta^*}
\Sigma^*_{\scriptscriptstyle\gamma\kappa}\,
\left(\frac{\partial \tau^{\scriptscriptstyle\beta}}{\partial\theta^{\scriptscriptstyle\,\kappa}}\right)_{\!\scriptscriptstyle \theta^*} \,=\,
\left(\frac{\partial \tau^{\scriptscriptstyle\alpha}}{\partial\theta^{\scriptscriptstyle\,\gamma}}\right)_{\!\scriptscriptstyle \theta^*}
\left(\frac{\partial \tau^{\scriptscriptstyle\beta}}{\partial\theta^{\scriptscriptstyle\,\gamma}}\right)_{\!\scriptscriptstyle \theta^*} \\[0.1cm]
\end{equation}
where the second equality follows because $\Sigma^*_{\scriptscriptstyle\gamma\kappa} = \delta_{\scriptscriptstyle\gamma\kappa}$ since $\Sigma^*$ is the identity matrix. 

It remains to show that $\Sigma^{\scriptscriptstyle\tau}$ is the inverse of the information matrix $I^{\scriptscriptstyle\tau}$ as in (\ref{eq:efficiency}). According to (\ref{eq:raofish2}), this is given by
\begin{equation} \label{eq:proffeff3}
I^{\scriptscriptstyle\tau}_{\scriptscriptstyle\alpha\beta} \,=\, 
\left.\frac{\partial^{\scriptscriptstyle\, 2}\!\,D}{\mathstrut\partial\tau^{\scriptscriptstyle \alpha}\partial\tau^{\scriptscriptstyle \beta}}\right|_{\scriptscriptstyle \tau^{\scriptscriptstyle \alpha} = 0} \,=\,
- \,E_{\scriptscriptstyle \theta^*}\left.\frac{\partial^{\scriptscriptstyle\, 2}\!\,\ell}{\mathstrut\partial\tau^{\scriptscriptstyle \alpha}\partial\tau^{\scriptscriptstyle \beta}}\right|_{\scriptscriptstyle \tau^{\scriptscriptstyle \alpha} = 0} \,=\,
E_{\scriptscriptstyle \theta^*}\! \left[\frac{\partial \ell}{\partial \tau^{\scriptscriptstyle\alpha}}\frac{\partial \ell}{\partial \tau^{\scriptscriptstyle\beta}}\right]_{\scriptscriptstyle \tau^{\scriptscriptstyle \alpha} = 0} \\[0.1cm]
\end{equation}
where the second equality follows from (\ref{eq:kl}), and the third equality from Fisher's identity (see~\cite{amari}, Page 28). Now, a direct application of the chain rule yields the following
$$
I^{\scriptscriptstyle\tau}_{\scriptscriptstyle\alpha\beta} \,=\, 
E_{\scriptscriptstyle \theta^*}\! \left[\frac{\partial \ell}{\partial \tau^{\scriptscriptstyle\alpha}}\frac{\partial \ell}{\partial \tau^{\scriptscriptstyle\beta}}\right]_{\scriptscriptstyle \tau^{\scriptscriptstyle \alpha} = 0} \,=\, 
\left(\frac{\partial \theta^{\scriptscriptstyle\,\gamma}}{\partial\tau^{\scriptscriptstyle\alpha}}\right)_{\!\scriptscriptstyle \theta^*}
E_{\scriptscriptstyle \theta^*}\! \left[\frac{\partial \ell}{\partial \theta^{\scriptscriptstyle\,\gamma}}\frac{\partial \ell}{\partial \theta^{\scriptscriptstyle\,\kappa}}\right]_{\scriptscriptstyle \theta^{\scriptscriptstyle\,\gamma} = 0}
\left(\frac{\partial \theta^{\scriptscriptstyle\,\kappa}}{\partial\tau^{\scriptscriptstyle\beta}}\right)_{\!\scriptscriptstyle \theta^*} \\[0.1cm]
$$
By the first equality in (\ref{eq:score2}), this is equal to
\begin{equation} \label{eq:proofeff4}
I^{\scriptscriptstyle\tau}_{\scriptscriptstyle\alpha\beta} \,=\, 
\left(\frac{\partial \theta^{\scriptscriptstyle\,\gamma}}{\partial\tau^{\scriptscriptstyle\alpha}}\right)_{\!\scriptscriptstyle \theta^*}
\Sigma^*_{\scriptscriptstyle\gamma\kappa}
\left(\frac{\partial \theta^{\scriptscriptstyle\,\kappa}}{\partial\tau^{\scriptscriptstyle\beta}}\right)_{\!\scriptscriptstyle \theta^*} \,=\,
\left(\frac{\partial \theta^{\scriptscriptstyle\,\gamma}}{\partial\tau^{\scriptscriptstyle\alpha}}\right)_{\!\scriptscriptstyle \theta^*}
\left(\frac{\partial \theta^{\scriptscriptstyle\,\gamma}}{\partial\tau^{\scriptscriptstyle\beta}}\right)_{\!\scriptscriptstyle \theta^*} \\[0.1cm]
\end{equation}
because $\Sigma^*_{\scriptscriptstyle\gamma\kappa} = \delta_{\scriptscriptstyle\gamma\kappa}$ is the identity matrix. Comparing (\ref{eq:proofeff2}) to (\ref{eq:proofeff4}), it is clear that $\Sigma^{\scriptscriptstyle\tau}$ is the inverse of the information matrix $I^{\scriptscriptstyle\tau}$ as in (\ref{eq:efficiency}).
\end{subequations}
\hfill$\blacksquare$ \\[0.1cm]
\textit{Proof of \emph{(iv)}}\,: (\ref{eq:information1}) and (\ref{eq:information2}) follow from (\ref{eq:ratel2}) and (\ref{eq:rateas}), respectively, by using (\ref{eq:raofish1}). Precisely, it is possible to write (\ref{eq:raofish1}) in the form
\begin{subequations}
\begin{equation} \label{eq:proofD1}
  D(\theta_{\scriptscriptstyle n}) \,=\, \frac{1}{2}\,d^{\scriptscriptstyle\,2}(\theta_{\scriptscriptstyle n},\theta^*) \,+\, \omega\!\left( \theta_{\scriptscriptstyle n}\right)d^{\scriptscriptstyle\,2}(\theta_{\scriptscriptstyle n},\theta^*)
\end{equation}
where $\omega$ is a continuous function in the neighborhood of $\theta^*$, equal to zero at $\theta = \theta^*$. To obtain (\ref{eq:information1}), it is enough to take expectations in (\ref{eq:proofD1}) and note that $\omega$ is bounded above in the neighborhood of $\theta^*$. Then, (\ref{eq:information1}) follows directly from (\ref{eq:ratel2}). 

To obtain (\ref{eq:information2}), it is enough to multiply (\ref{eq:proofD1}) by $n^{\scriptscriptstyle p}$ where $p \in (0,1)$. This gives the following expression
\begin{equation} \label{eq:proofD2}
n^{\scriptscriptstyle p}  D(\theta_{\scriptscriptstyle n}) \,=\,
\frac{1}{2}\,n^{\scriptscriptstyle p}d^{\scriptscriptstyle\,2}(\theta_{\scriptscriptstyle n},\theta^*) \left(1+
\omega\!\left( \theta_{\scriptscriptstyle n}\right)\right)
\end{equation}
From (\ref{eq:rateas}), $n^{\scriptscriptstyle p}d^{\scriptscriptstyle\,2}(\theta_{\scriptscriptstyle n},\theta^*)$ converges to zero almost surely. Moreover, by continuity of $\omega$, it follows that $\omega\!\left( \theta_{\scriptscriptstyle n}\right)$ converges to $\omega\!\left( \theta^*\right) = 0$ almost surely. Therefore, by taking limits in (\ref{eq:proofD2}), it is readily seen that 
\begin{equation}
\lim\, n^{\scriptscriptstyle p}  D(\theta_{\scriptscriptstyle n}) \,=\,
\frac{1}{2}\left(\lim\, n^{\scriptscriptstyle p}d^{\scriptscriptstyle\,2}(\theta_{\scriptscriptstyle n},\theta^*)\right)
\left(1+
\lim\,\omega\!\left( \theta_{\scriptscriptstyle n}\right) \right) \,=\, 0 
\end{equation}
almost surely. However, this is equivalent to the statement that $D(\theta_{\scriptscriptstyle n}) = o(n^{\scriptscriptstyle -p})$ for $p \in (0,1)$, almost surely. Thus,
(\ref{eq:information2}) is proved. \hfill$\blacksquare$ 
\end{subequations}
\vfill
\pagebreak
\appendix

\section{Proofs of geometric lemmas} \label{sec:geometric} 
\subsection{Lemma \ref{lemma:grad}} let $c(t)$ be the geodesic connecting $\theta^*$ to some $\theta \in \Theta^*$, parameterised by arc length. In other words, $c(0) = \theta^*$ and $c(t_{\scriptscriptstyle\theta}) = \theta$ where $t_{\scriptscriptstyle\theta} = d(\theta,\theta^*)$. Denote $\Pi_{\scriptscriptstyle t}$ the parallel transport along $c(t)$, from $T_{\scriptscriptstyle c(0)}\Theta$ to $T_{\scriptscriptstyle c(t)}\Theta$. Since the velocity $\dot{c}(t)$ is self-parallel~\cite{petersen},
$$
\dot{c}(t_{\scriptscriptstyle\theta}) \,=\, \Pi_{\scriptscriptstyle t_{\scriptscriptstyle\theta}}(\dot{c}(0))
$$
Multiplying this identity by $-t_{\scriptscriptstyle\theta\,}$, it follows that
\begin{subequations} 
\begin{equation} \label{eq:proofgrad1}
\mathrm{Exp}^{\scriptscriptstyle -1}_{\scriptscriptstyle \theta}(\theta^*) \,=\, 
-\,\Pi_{\scriptscriptstyle t_{\scriptscriptstyle\theta}}\!\left(\mathrm{Exp}^{\scriptscriptstyle -1}_{\scriptscriptstyle \theta^*}(\theta)\right)
\end{equation}
Moreover, recall the first-order Taylor expansion of the gradient $\nabla D(\theta)$~\cite{petersen}\cite{chavel}
\begin{equation} \label{eq:proofgrad2}
   \nabla D(\theta) \,=\, \Pi_{\scriptscriptstyle t_{\scriptscriptstyle\theta}}\!\left(\nabla D(\theta^*) + t_{\scriptscriptstyle\theta}\,\nabla^2D(\theta^*)\cdot\dot{c}(0)+t_{\scriptscriptstyle\theta}\,\phi(\theta)\right)
\end{equation}
where $\phi(\theta)$ is continuous and equal to zero at $\theta = \theta^*$. Here, $\nabla^2D(\theta^*)$ is the Hessian of $D(\theta)$ at $\theta = \theta^*$, considered as a linear mapping of $T_{\scriptscriptstyle\theta^*}\Theta$~\cite{petersen}\cite{chavel}
$$
\nabla^2D(\theta^*)\cdot w \,=\, \nabla_{\scriptscriptstyle w} \nabla D(\theta^*) \hspace{1cm} \text{for }\,\,w \in T_{\scriptscriptstyle\theta^*}\Theta
$$
where $\nabla_{\scriptscriptstyle w}$ denotes the covariant derivative in the direction of $w$. By (d1), the first term on the right-hand side of (\ref{eq:proofgrad2}) is equal to zero, so that
\begin{equation} \label{eq:proofgrad3}
\nabla D(\theta) \,=\, \Pi_{\scriptscriptstyle t_{\scriptscriptstyle\theta}}\!\left(\nabla^2D(\theta^*)\cdot\mathrm{Exp}^{\scriptscriptstyle -1}_{\scriptscriptstyle \theta^*}(\theta)+t_{\scriptscriptstyle\theta}\,\phi(\theta)\right)
\end{equation} 
Taking the scalar product of (\ref{eq:proofgrad1}) and (\ref{eq:proofgrad3}),
\begin{equation} \label{eq:proofgrad4}
  \langle\mathrm{Exp}^{\scriptscriptstyle -1}_{\scriptscriptstyle \theta}(\theta^*),\nabla D(\theta)\rangle\,=\,
  -\,\langle\mathrm{Exp}^{\scriptscriptstyle -1}_{\scriptscriptstyle \theta^*}(\theta),\nabla^2D(\theta^*)\cdot
\mathrm{Exp}^{\scriptscriptstyle -1}_{\scriptscriptstyle \theta^*}(\theta)\rangle \,-\, t_{\scriptscriptstyle\theta}\,
\langle\mathrm{Exp}^{\scriptscriptstyle -1}_{\scriptscriptstyle \theta^*}(\theta),\phi(\theta)\rangle
\end{equation}
since parallel transport preserves scalar products. In terms of the normal coordinates $\theta^{\scriptscriptstyle\,\alpha}$, this reads~\cite{petersen}
\begin{equation} \label{eq:proofgrad5}
\langle\mathrm{Exp}^{\scriptscriptstyle -1}_{\scriptscriptstyle \theta}(\theta^*),\nabla D(\theta)\rangle\,=\,
-\,H_{\scriptscriptstyle\alpha\beta}\,\theta^{\scriptscriptstyle\,\alpha}\theta^{\scriptscriptstyle\,\beta}\,-\,
t^{\scriptscriptstyle2}_{\scriptscriptstyle\theta}\,\,\hat{\theta}^{\scriptscriptstyle\,\alpha}\phi^{\scriptscriptstyle\alpha}
\end{equation}
where $H = (H_{\scriptscriptstyle\alpha\beta})$ was defined in (\ref{eq:hess}), $\hat{\theta}^{\scriptscriptstyle\,\alpha}$ denotes the quotient $\left.\theta^{\scriptscriptstyle\,\alpha}\middle/t_{\scriptscriptstyle\theta\,}\right.$, and the $\phi^{\scriptscriptstyle\alpha}$ denote the components of $\phi(\theta)$. Note that $t^{\scriptscriptstyle2}_{\scriptscriptstyle\theta} = d^{\scriptscriptstyle\,2}(\theta,\theta^*) = \theta^{\scriptscriptstyle\,\alpha}\theta^{\scriptscriptstyle\,\alpha\,}$, so (\ref{eq:proofgrad5}) can be written
\begin{equation} \label{eq:proofgrad6}
\langle\mathrm{Exp}^{\scriptscriptstyle -1}_{\scriptscriptstyle \theta}(\theta^*),\nabla D(\theta)\rangle\,=\,
\left( \psi(\theta)\delta_{\scriptscriptstyle\alpha\beta} \,-\,H_{\scriptscriptstyle\alpha\beta}\right)\,\theta^{\scriptscriptstyle\,\alpha}\theta^{\scriptscriptstyle\,\beta}
\end{equation}
where $\psi(\theta)$ is continuous and equal to zero at $\theta = \theta^*$. To conclude, let $\mu = \lambda - \varepsilon$ for some $\varepsilon > 0$, and $\bar{\Theta}^*$ a neighborhood of $\theta^*$, contained in $\Theta^*$, such that $\psi(\theta) \leq \varepsilon$ for $\theta \in \bar{\Theta}^*$. Then, since $\lambda$ is the smallest eigenvalue of $H = (H_{\scriptscriptstyle\alpha\beta})$,
$$
\langle\mathrm{Exp}^{\scriptscriptstyle -1}_{\scriptscriptstyle \theta}(\theta^*),\nabla D(\theta)\rangle\,\leq\,
\left( \varepsilon - \lambda\right)\,\theta^{\scriptscriptstyle\,\alpha}\theta^{\scriptscriptstyle\,\alpha} \,=\, -\mu\,d^{\scriptscriptstyle\,2}(\theta,\theta^*)
$$
for $\theta \in \bar{\Theta}^*$. This is exactly (\ref{eq:lemmgrad}), so the lemma is proved.\hfill$\blacksquare$
\end{subequations}
\subsection{Lemma \ref{lemma:linearalgo}} to simplify notation, let $u_{\scriptscriptstyle n+1}=u(\theta_{\scriptscriptstyle n},x_{\scriptscriptstyle n+1})$. Then, the geodesic $c(t)$, connecting $\theta_{\scriptscriptstyle n}$ to $\theta_{\scriptscriptstyle n+1\,}$, has equation
$$
c(t) \,=\, \mathrm{Exp}_{\scriptscriptstyle \theta_{\scriptscriptstyle n}}\!\left(t\gamma_{\scriptscriptstyle n+1}u_{\scriptscriptstyle n+1}\right)
$$
Each one of the normal coordinates $\theta^{\scriptscriptstyle\,\alpha}$ is a $C^3$ function $\theta^{\scriptscriptstyle\,\alpha}:\Theta^*\rightarrow \mathbb{R}$, with differential $d\theta^{\scriptscriptstyle\,\alpha}$ and Hessian~\cite{petersen}
$$
\nabla^2\theta^{\scriptscriptstyle\,\alpha}\,=\, - \Gamma^{\scriptscriptstyle\alpha}_{\scriptscriptstyle\beta\gamma}(\theta)\,d\theta^{\scriptscriptstyle\,\beta\,}{\scriptstyle\otimes}\,\,d\theta^{\scriptscriptstyle\,\gamma}
$$
where $\Gamma^{\scriptscriptstyle\alpha}_{\scriptscriptstyle\beta\gamma}$ are the Christoffel symbols of the coordinates $\theta^{\scriptscriptstyle\,\alpha}$, and ${\scriptstyle\otimes}$ denotes the tensor product. Then, the second-order Taylor expansion of the functions $\theta^{\scriptscriptstyle\,\alpha}\circ c$ reads
\begin{subequations}
\begin{equation} \label{eq:prooflem31}
   (\theta^{\scriptscriptstyle\,\alpha}\circ c)(1) =    (\theta^{\scriptscriptstyle\,\alpha}\circ c)(0)
+\gamma_{\scriptscriptstyle n+1}\,d\theta^{\scriptscriptstyle\,\alpha}(u_{\scriptscriptstyle n+1})
-\frac{1}{2}\gamma^{\scriptscriptstyle 2}_{\scriptscriptstyle n+1}\,\Gamma^{\scriptscriptstyle\alpha}_{\scriptscriptstyle\beta\gamma}(\theta_{\scriptscriptstyle n})\,d\theta^{\scriptscriptstyle\,\beta}(u_{\scriptscriptstyle n+1})\,d\theta^{\scriptscriptstyle\,\gamma}(u_{\scriptscriptstyle n+1})
+\,\gamma^{\scriptscriptstyle 3}_{\scriptscriptstyle n+1}T^{\scriptscriptstyle\,\alpha}_{\scriptscriptstyle n+1}
\end{equation}
where $T^{\scriptscriptstyle\,\alpha}_{\scriptscriptstyle n+1}$ satisfies
\begin{equation}  \label{eq:prooflem32}
  \left|T^{\scriptscriptstyle\,\alpha}_{\scriptscriptstyle n+1}\right| \,\leq\,K_{\scriptscriptstyle 1}\,\Vert u_{\scriptscriptstyle n+1}\Vert^3
\end{equation}
for a constant $K_{\scriptscriptstyle 1}$ which does not depend on $n$, as can be shown by direct calculation. Of course, $(\theta^{\scriptscriptstyle\,\alpha}\circ c)(1) = \theta^{\scriptscriptstyle\,\alpha}_{\scriptscriptstyle n+1}$ and $(\theta^{\scriptscriptstyle\,\alpha}\circ c)(0) = \theta^{\scriptscriptstyle\,\alpha}_{\scriptscriptstyle n\,}$. Moreover, $d\theta^{\scriptscriptstyle\,\alpha}(u_{\scriptscriptstyle n+1}) = u^{\scriptscriptstyle\alpha}_{\scriptscriptstyle n+1}$ are the components of $u_{\scriptscriptstyle n+1\,}$. Replacing into (\ref{eq:prooflem31}), this yields
\begin{equation} \label{eq:prooflem33}
\theta^{\scriptscriptstyle\,\alpha}_{\scriptscriptstyle n+1} \,=\, \theta^{\scriptscriptstyle\,\alpha}_{\scriptscriptstyle n}\,+\,\gamma^{\phantom{\scriptscriptstyle 2}}_{\scriptscriptstyle n+1}\,u^{\scriptscriptstyle\alpha}_{\scriptscriptstyle n+1}\,+\,
\gamma^{\scriptscriptstyle 2}_{\scriptscriptstyle n+1}\,\pi^{\scriptscriptstyle\alpha}_{\scriptscriptstyle n+1}
\end{equation}
where $\pi^{\scriptscriptstyle\alpha}_{\scriptscriptstyle n+1}$ is given by
\begin{equation} \label{eq:prooflem34}
\pi^{\scriptscriptstyle\alpha}_{\scriptscriptstyle n+1}\,=\, 
\gamma^{\phantom{\scriptscriptstyle 2}}_{\scriptscriptstyle n+1}\,T^{\scriptscriptstyle\,\alpha}_{\scriptscriptstyle n+1}\,-\,\frac{1}{2}\,\Gamma^{\scriptscriptstyle\alpha}_{\scriptscriptstyle\beta\gamma}(\theta_{\scriptscriptstyle n})\,u^{\scriptscriptstyle\beta}_{\scriptscriptstyle n+1}\,u^{\scriptscriptstyle\gamma}_{\scriptscriptstyle n+1}
\end{equation}
Comparing (\ref{eq:prooflem33}) to (\ref{eq:linearalgo}), it is clear the proof will be complete upon showing $\mathbb{E}\left|\pi^{\scriptscriptstyle \alpha}_{\scriptscriptstyle n+1}\right| = O(n^{\scriptscriptstyle -1/2})$. To do so,  note that each Christoffel symbol $\Gamma^{\scriptscriptstyle\alpha}_{\scriptscriptstyle\beta\gamma}$ is a $C^1$ function on the compact set $\Theta^*$, with $\Gamma^{\scriptscriptstyle\alpha}_{\scriptscriptstyle\beta\gamma}(\theta^*) = 0$ by the definition of normal coordinates~\cite{petersen}. Therefore,
\begin{equation} \label{eq:pi1}
  \left|\,\Gamma^{\scriptscriptstyle\alpha}_{\scriptscriptstyle\beta\gamma}(\theta)\,\right| \,\leq\,K_{\scriptscriptstyle 2}\,d(\theta,\theta^*)
\end{equation}
for a constant $K_{\scriptscriptstyle 2}$ which does not depend on $n$. Replacing the inequalities (\ref{eq:prooflem32}) and (\ref{eq:pi1}) into (\ref{eq:prooflem34}), and taking expectations, it follows that
\end{subequations}
\begin{subequations}
  \begin{equation} \label{eq:ep1}
\mathbb{E}\left|\pi^{\scriptscriptstyle \alpha}_{\scriptscriptstyle n+1}\right|\,\leq\,
\gamma_{\scriptscriptstyle n+1}\,K_{\scriptscriptstyle 1}\,\mathbb{E}\,\Vert u_{\scriptscriptstyle n+1}\Vert^3\,+\,d^{\scriptscriptstyle\,2}\times K_{\scriptscriptstyle 2}\,\mathbb{E}\left[\,d(\theta_{\scriptscriptstyle n},\theta^*)\,
\Vert u_{\scriptscriptstyle n+1}\Vert^2\,\right]
  \end{equation}
where $d$ is the dimension of the parameter space $\Theta$. However, using the fact that the $x_{\scriptscriptstyle n}$ are i.i.d. with distribution $P_{\scriptscriptstyle \theta^*\,}$,
\begin{equation} \label{eq:ep2}
   \mathbb{E}\left[\,\Vert u_{\scriptscriptstyle n+1}\Vert^3\,\middle|\mathcal{X}_{\scriptscriptstyle n}\,\right]\,=\,E_{\scriptscriptstyle\theta^*}\,\Vert u(\theta_{\scriptscriptstyle n},x)\Vert^3\,\leq\, R^{\scriptscriptstyle\, 3/4}(\theta_{\scriptscriptstyle n})
\end{equation}
by (u2) and Jensen's inequality~\cite{shiryayev}. On the other hand, by the Cauchy-Schwarz inequality,
$$
\mathbb{E}\left[\,d(\theta_{\scriptscriptstyle n},\theta^*)\,
\Vert u_{\scriptscriptstyle n+1}\Vert^2\,\right] \,\leq\, \left(\mathbb{E}\,d^{\scriptscriptstyle\,2}(\theta_{\scriptscriptstyle n},\theta^*)\right)^{\scriptscriptstyle 1/2}\,\left(\mathbb{E}\,\Vert u_{\scriptscriptstyle n+1}\Vert^{\scriptscriptstyle 4}\right)^{\scriptscriptstyle 1/2}
\,\leq\,b\,n^{\scriptscriptstyle-1/2}\,\left(\mathbb{E}\,\Vert u_{\scriptscriptstyle n+1}\Vert^{\scriptscriptstyle 4}\right)^{\scriptscriptstyle 1/2}
$$
for some $b > 0$ as follows from (\ref{eq:ratel2}). Then, by the same reasoning that lead to (\ref{eq:ep2}),
\begin{equation} \label{eq:ep3}
\mathbb{E}\left[\,d(\theta_{\scriptscriptstyle n},\theta^*)\,
\Vert u_{\scriptscriptstyle n+1}\Vert^2\,\right] \,\leq\, b\,n^{\scriptscriptstyle-1/2}\,\left(\mathbb{E}\,R(\theta_{\scriptscriptstyle n})\right)^{\scriptscriptstyle 1/2}
\end{equation}
By (u2), there exists a uniform upper bound $M$ on $R(\theta)$ for $\theta \in \Theta^*$. Since $\theta_{\scriptscriptstyle n}$ lies in $\Theta^*$ for all $n$, it follows by replacing the inequalities (\ref{eq:ep2}) and (\ref{eq:ep3}) into (\ref{eq:ep1}) that
\begin{equation} \label{eq:ep4}
 \mathbb{E}\left|\pi^{\scriptscriptstyle \alpha}_{\scriptscriptstyle n+1}\right|\,\leq\,
\gamma_{\scriptscriptstyle n+1}\,K_{\scriptscriptstyle 1}\,M^{\scriptscriptstyle\,3/4}\,\,+\,
\,d^{\scriptscriptstyle\,2}\times K_{\scriptscriptstyle 2}\,b\,n^{\scriptscriptstyle-1/2}M^{\scriptscriptstyle\,1/2}
\end{equation}
Finally, by recalling that $\gamma_{\scriptscriptstyle n} = \frac{a}{n\,}$, it is clear that the right-hand side of (\ref{eq:ep4}) is $O(n^{\scriptscriptstyle -1/2})$, so the proof is complete.\hfill$\blacksquare$
\end{subequations}
\subsection{Lemma \ref{lemma:linearfield}} the lemma is an instance of the general statement\,: let $\theta \in \Theta^*$ and $v = \nabla D(\theta)$. Then, in a system of normal coordinates with origin at $\theta^*$, 
\begin{subequations}
\begin{equation}\label{eq:prooflf1}
  v^{\scriptscriptstyle\,\alpha}\,=\, H_{\scriptscriptstyle\alpha\beta}\,\theta^{\scriptscriptstyle\,\beta}\,+\, o\left(d(\theta,\theta^*)\right)
\end{equation}
where $v^{\scriptscriptstyle\,\alpha}$ are the components of $v$. Indeed, (\ref{eq:linearfield}) follows from (\ref{eq:prooflf1}) after replacing $\theta = \theta_{\scriptscriptstyle n\,}$, so that $v = v_{\scriptscriptstyle n\,}$, and setting
$$
\rho^{\scriptscriptstyle\alpha}_{\scriptscriptstyle n}\,=\, v^{\scriptscriptstyle\,\alpha}_{\scriptscriptstyle n}\,-\,H_{\scriptscriptstyle\alpha\beta}\,\theta^{\scriptscriptstyle\, \beta}_{\scriptscriptstyle n}
$$
To prove (\ref{eq:prooflf1}), recall (\ref{eq:proofgrad3}) from the proof of Lemma \ref{lemma:grad}, which can be written
\begin{equation} \label{eq:prooflf2}
  v\,=\, 
\Pi_{\scriptscriptstyle t_{\scriptscriptstyle\theta}}\!\left(\nabla^2D(\theta^*)\cdot\mathrm{Exp}^{\scriptscriptstyle -1}_{\scriptscriptstyle \theta^*}(\theta)\right)\,+\,
d(\theta,\theta^*)\,\Pi_{\scriptscriptstyle t_{\scriptscriptstyle\theta}}\!\left(\phi(\theta)\right)
\end{equation}
Denote $\partial_{\scriptscriptstyle\alpha}=\frac{\partial}{\mathstrut\partial\theta^{\scriptscriptstyle\,\alpha}}$ the coordinate vector fields of the normal coordinates $\theta^{\scriptscriptstyle\,\alpha\,}$. Note that~\cite{petersen}\cite{chavel}
$$
\mathrm{Exp}^{\scriptscriptstyle -1}_{\scriptscriptstyle \theta^*}(\theta)\,=\,\theta^{\scriptscriptstyle\,\beta}\,
\partial_{\scriptscriptstyle\beta}(\theta^*) \hspace{1cm}
\nabla^2D(\theta^*)\cdot \partial_{\scriptscriptstyle\beta}(\theta^*) \,=\,H_{\scriptscriptstyle\alpha\beta}\,\partial_{\scriptscriptstyle\alpha}(\theta^*)
$$
Replacing in (\ref{eq:prooflf2}), this gives
\begin{equation} \label{eq:prooflf3}
v\,=\,H_{\scriptscriptstyle\alpha\beta}\,\theta^{\scriptscriptstyle\,\beta}\,\,\Pi_{\scriptscriptstyle t_{\scriptscriptstyle\theta}}\!\left(\partial_{\scriptscriptstyle\alpha}(\theta^*)\right)\,+\,d(\theta,\theta^*)\,\Pi_{\scriptscriptstyle t_{\scriptscriptstyle\theta}}\!\left(\phi(\theta)\right)
\end{equation}
From the first-order Taylor expansion of the vector fields $\partial_{\scriptscriptstyle\alpha}$~\cite{petersen}\cite{chavel}
$$
 \partial_{\scriptscriptstyle\alpha}(\theta)\,=\,\Pi_{\scriptscriptstyle t_{\scriptscriptstyle\theta}}\!\left(\partial_{\scriptscriptstyle\alpha}(\theta^*)\,+\,\nabla \partial_{\scriptscriptstyle\alpha}(\theta^*)\cdot\mathrm{Exp}^{\scriptscriptstyle -1}_{\scriptscriptstyle \theta^*}(\theta) \right)\,+\,d(\theta,\theta^*)\,\Pi_{\scriptscriptstyle t_{\scriptscriptstyle\theta}}\!\left(\chi^{\scriptscriptstyle\alpha}(\theta)\right)
$$
where $\chi^{\scriptscriptstyle\alpha}(\theta)$ is continuous and equal to zero at $\theta = \theta^*$. However, by the definition of normal coordinates~\cite{petersen}, each covariant derivative $\nabla \partial_{\scriptscriptstyle\alpha}(\theta^*)$ is zero. In other words,
\begin{equation} \label{eq:prooflf4}
 \partial_{\scriptscriptstyle\alpha}(\theta)\,=\,\Pi_{\scriptscriptstyle t_{\scriptscriptstyle\theta}}\!\left(\partial_{\scriptscriptstyle\alpha}(\theta^*)\right)\,+\,d(\theta,\theta^*)\,\Pi_{\scriptscriptstyle t_{\scriptscriptstyle\theta}}\!\left(\chi^{\scriptscriptstyle\alpha}(\theta)\right)
\end{equation}
Replacing (\ref{eq:prooflf4}) into (\ref{eq:prooflf3}), it follows
\begin{equation} 
v\,=\,H_{\scriptscriptstyle\alpha\beta}\,\theta^{\scriptscriptstyle\,\beta}\,\partial_{\scriptscriptstyle\alpha}(\theta)\,+\,d(\theta,\theta^*)\,
\Pi_{\scriptscriptstyle t_{\scriptscriptstyle\theta}}\!\left(\phi(\theta) - H_{\scriptscriptstyle\alpha\beta}\,\theta^{\scriptscriptstyle\,\beta}\chi^{\scriptscriptstyle\alpha}(\theta)\right)
\end{equation}
\end{subequations}
Now, to obtain (\ref{eq:prooflf1}), it is enough to note the decomposition $v = v^{\scriptscriptstyle\,\alpha}\,\partial_{\scriptscriptstyle\alpha}(\theta)$ is unique, and $\phi(\theta) - H_{\scriptscriptstyle\alpha\beta}\,\theta^{\scriptscriptstyle\,\beta}\chi^{\scriptscriptstyle\alpha}(\theta)$ converges to zero as $\theta$ converges to $\theta^*$. \hfill$\blacksquare$

\section{Conditions of the martingale CLT} \label{sec:clt}
for the verification of Conditions (\ref{subeq:clt}), the following inequality (\ref{eq:stableA}) will be useful. Let $\nu = a\lambda - \frac{1}{2\,}$, so $-\nu$ is the largest eigenvalue of the matrix $A$ in (\ref{eq:linearisation3}). There exists a constant $C_{\scriptscriptstyle\! A}$ such that the transition matrices $A_{\scriptscriptstyle n,k}$ in (\ref{eq:transitions2}) satisfy~\cite{nev}\cite{kailath}
\begin{equation} \label{eq:stableA}
  \left| A_{\scriptscriptstyle n,k} \right|_{\scriptscriptstyle \mathrm{Op}}\,\leq C_{\scriptscriptstyle\! A}\left(\frac{k}{n}\right)^{\!\nu}
\end{equation}
where $\left|A_{\scriptscriptstyle n,k}\right|_{\scriptscriptstyle \mathrm{Op}}$ denotes the Euclidean operator norm, equal to the largest singular value of the matrix $A_{\scriptscriptstyle n,k\,}$. \\[0.1cm]
\textit{Condition (\ref{eq:clt1})}\,: to verify this condition, note that for arbitrary $\varepsilon > 0$,
\begin{subequations}
\begin{equation} \label{eq:clt11}
\mathbb{P}\left(
\max_{\scriptscriptstyle k\leq n}\,\left| A^{\phantom{\scriptscriptstyle\,2}}_{\scriptscriptstyle n,k}\,\frac{aw_{\scriptscriptstyle k}}{k^{\scriptscriptstyle 1/2}}\right| > \varepsilon \right) \,\leq\, \sum^{\scriptscriptstyle n}_{\scriptscriptstyle k=1}\,\mathbb{P}
\left(\left| A^{\phantom{\scriptscriptstyle\,2}}_{\scriptscriptstyle n,k}\,\frac{aw_{\scriptscriptstyle k}}{k^{\scriptscriptstyle 1/2}}\right| > \varepsilon \right) \,\leq\, 
\sum^{\scriptscriptstyle n}_{\scriptscriptstyle k=1}\,\mathbb{P}
\left(C_{\scriptscriptstyle\! A}\left(\frac{k}{n}\right)^{\!\!\nu}\left|\frac{aw_{\scriptscriptstyle k}}{k^{\scriptscriptstyle 1/2}}\right| > \varepsilon \right)
\end{equation}
where the second inequality follows from (\ref{eq:stableA}). However, it follows from (u2) that there exists a uniform upper bound $M_{\scriptscriptstyle w}$ on the fourth-order moments of $|w_{\scriptscriptstyle k}|\,$. Therefore, by Chebyshev's inequality~\cite{shiryayev}
\begin{equation} \label{eq:clt12}
\sum^{\scriptscriptstyle n}_{\scriptscriptstyle k=1}\,\mathbb{P}
\left(C_{\scriptscriptstyle\! A}\left(\frac{k}{n}\right)^{\!\!\nu}\left|\frac{aw_{\scriptscriptstyle k}}{k^{\scriptscriptstyle 1/2}}\right| > \varepsilon \right) \,\leq\, \left(\frac{aC_{\scriptscriptstyle\! A}}{\varepsilon}\right)^{\!4}\frac{M_{\scriptscriptstyle w}}{n^{\scriptscriptstyle 4\nu}}\,\sum^{\scriptscriptstyle n}_{\scriptscriptstyle k=1}\,k^{\scriptscriptstyle 4\nu-2}
\end{equation}
\end{subequations}
Since $\nu > 0$, the right-hand side of (\ref{eq:clt12}) has limit equal to $0$ as $n \rightarrow \infty$, by the Euler-Maclaurin formula~\cite{courant}. Replacing this limit from (\ref{eq:clt12}) into (\ref{eq:clt11}) immediately yields Condition (\ref{eq:clt1}). \hfill$\blacksquare$ \\[0.1cm]
\textit{Condition (\ref{eq:clt2})}\,: to verify this condition, recall that $(w_{\scriptscriptstyle k})$ is a sequence of square-integrable martingale differences. Therefore, from (\ref{eq:sum})
\begin{subequations}
\begin{equation} \label{eq:clt21}
   \mathbb{E}\left|\tilde{\eta}^{\phantom{\scriptscriptstyle\,2}}_{\scriptscriptstyle n}\right|^2\,=\,\sum^{\scriptscriptstyle n}_{\scriptscriptstyle k=1}\,\frac{a^{\scriptscriptstyle 2}}{k}\,\mathbb{E}\,\mathrm{tr}\!\left(A^{\scriptscriptstyle 2}_{\scriptscriptstyle n,k}\Sigma^{\phantom{\scriptscriptstyle 2}}_{\scriptscriptstyle k} \right)
\end{equation}
where $\Sigma_{\scriptscriptstyle k}$ is the conditional covariance matrix in (\ref{eq:sigmak}). Applying (\ref{eq:stableA}) to each term under the sum in (\ref{eq:clt21}), it follows that
\begin{equation} \label{eq:clt22}
\mathbb{E}\left|\tilde{\eta}^{\phantom{\scriptscriptstyle\,2}}_{\scriptscriptstyle n}\right|^2\,\leq\,
d^{\scriptscriptstyle\frac{1}{\mathstrut 2}}\,\sum^{\scriptscriptstyle n}_{\scriptscriptstyle k=1}\,\frac{a^{\scriptscriptstyle 2}}{k}\,\mathbb{E}\,\left| A_{\scriptscriptstyle n,k} \right|^{\scriptscriptstyle 2}_{\scriptscriptstyle \mathrm{Op}}\left|\Sigma^{\phantom{\scriptscriptstyle 2}}_{\scriptscriptstyle k}\right|_{\scriptscriptstyle\mathrm{F}}
\,\leq\,
\left(d^{\scriptscriptstyle\frac{1}{\mathstrut 2}}\,a^{\scriptscriptstyle 2}\,C^{\scriptscriptstyle\, 2\,}_{\scriptscriptstyle\! A}\right)\,\frac{1}{n^{\scriptscriptstyle 2\nu}}\,\sum^{\scriptscriptstyle n}_{\scriptscriptstyle k=1}\,k^{\scriptscriptstyle 2\nu-1}\,\mathbb{E}\left|\Sigma^{\phantom{\scriptscriptstyle 2}}_{\scriptscriptstyle k}\right|_{\scriptscriptstyle\mathrm{F}}
\end{equation}
where $d$ is the dimension of  the parameter space $\Theta$, and $\left|\Sigma^{\phantom{\scriptscriptstyle 2}}_{\scriptscriptstyle k}\right|_{\scriptscriptstyle\mathrm{F}}$ denotes the Frobenius matrix norm. However, it follows from (u1) that there exists a uniform upper bound $S$ on $\left|\Sigma^{\phantom{\scriptscriptstyle 2}}_{\scriptscriptstyle k}\right|_{\scriptscriptstyle\mathrm{F}\,}$. Therefore, by (\ref{eq:clt22})
\begin{equation} \label{eq:clt23}
 \mathbb{E}\left|\tilde{\eta}^{\phantom{\scriptscriptstyle\,2}}_{\scriptscriptstyle n}\right|^2\,\leq\, 
\left(d^{\scriptscriptstyle\frac{1}{\mathstrut 2}}\,a^{\scriptscriptstyle 2}\,C^{\scriptscriptstyle\, 2\,}_{\scriptscriptstyle\! A}\right)\,\frac{S}{n^{\scriptscriptstyle 2\nu}}\,\sum^{\scriptscriptstyle n}_{\scriptscriptstyle k=1}\,k^{\scriptscriptstyle 2\nu-1}
\end{equation}
Since $\nu > 0$, the right-hand side of (\ref{eq:clt23}) remains bounded as $n \rightarrow \infty$, by the Euler-Maclaurin formula~\cite{courant}. This immediately yields Condition (\ref{eq:clt2}).
\end{subequations} \hfill$\blacksquare$ \\[0.1cm]
\textit{Condition (\ref{eq:clt3})}\,: to verify this condition, it is first admitted that the following limit is known to hold
\begin{subequations}
\begin{equation} \label{eq:clt31}
  \lim\,\mathbb{E}\left(\Sigma_{\scriptscriptstyle k}\right) \,=\, \Sigma^*
\end{equation}
where $\Sigma^*$ was defined in (\ref{eq:cov}). Then, let the sum in (\ref{eq:clt3}) be written
\begin{equation} \label{eq:riemann1}
   \sum^{\scriptscriptstyle n}_{\scriptscriptstyle k=1}\frac{a^{\scriptscriptstyle 2}}{k}\,A^{\phantom{\scriptscriptstyle\,2}}_{\scriptscriptstyle n,k\,}\Sigma^{\phantom{\scriptscriptstyle\,2}}_{\scriptscriptstyle k\,}A^{\phantom{\scriptscriptstyle\,2}}_{\scriptscriptstyle n,k} \,=\,
\sum^{\scriptscriptstyle n}_{\scriptscriptstyle k=1}\frac{a^{\scriptscriptstyle 2}}{k}\,A^{\phantom{\scriptscriptstyle\,2}}_{\scriptscriptstyle n,k\,}\Sigma^*A^{\phantom{\scriptscriptstyle\,2}}_{\scriptscriptstyle n,k}\,+\,
\sum^{\scriptscriptstyle n}_{\scriptscriptstyle k=1}\frac{a^{\scriptscriptstyle 2}}{k}\,A^{\phantom{\scriptscriptstyle\,2}}_{\scriptscriptstyle n,k\,}\left[\,\Sigma^{\phantom{\scriptscriptstyle\,2}}_{\scriptscriptstyle k} - \Sigma^*\right]A^{\phantom{\scriptscriptstyle\,2}}_{\scriptscriptstyle n,k}
\end{equation}
Due to the equivalence $A_{\scriptscriptstyle n,k} \sim \exp(\ln(n/k)A)$ (see~\cite{nev}, Page 125),  the first term in (\ref{eq:riemann1}) is a Riemann sum for the integral~\cite{nev}\cite{kailath}
$$
a^{\scriptscriptstyle 2}\,\int^1_0\,e^{\scriptscriptstyle -\ln(s)\,A}\,\Sigma^*\,e^{\scriptscriptstyle -\ln(s)\,A}\,d\ln(s) \,=\,
a^{\scriptscriptstyle 2}\,\int^\infty_0\,e^{\scriptscriptstyle -t\,A}\,\Sigma^*\,e^{\scriptscriptstyle -t\,A}\,dt 
$$
which is known to be the solution $\Sigma$ of Lyapunov's equation (\ref{eq:lyapunov}). The second term in (\ref{eq:riemann1}) can be shown to converge to zero in probability, using inequality (\ref{eq:stableA}) and the limit (\ref{eq:clt31}), by a similar argument to the ones in the verification of Conditions (\ref{eq:clt1}) and (\ref{eq:clt2}). Then, Condition (\ref{eq:clt3}) follows immediately. \hfill$\blacksquare$ \\[0.1cm]
\end{subequations}
\textit{Proof of (\ref{eq:clt31})}\,:
\begin{subequations}
recall that $w_{\scriptscriptstyle k} = u_{\scriptscriptstyle k} + v_{\scriptscriptstyle k-1}$ where $u_{\scriptscriptstyle k} = u(\theta_{\scriptscriptstyle k-1},x_{\scriptscriptstyle k})$ and $v_{\scriptscriptstyle k-1} = \nabla D(\theta_{\scriptscriptstyle k-1})$. Since $(w_{\scriptscriptstyle k})$ is a sequence of square-integrable martingale differences, it is possible to write, in the notation of (\ref{eq:sigmak}),
\begin{equation} \label{eq:clt32}
  \Sigma^{\phantom{\scriptscriptstyle2}}_{\scriptscriptstyle k} \,=\, \mathbb{E}\left[\,u^{\phantom{\scriptscriptstyle{\dagger}}}_{\scriptscriptstyle k}u^{\scriptscriptstyle{\dagger}}_{\scriptscriptstyle k}\,\middle|\mathcal{X}_{\scriptscriptstyle k-1}\right]\,-\, v^{\phantom{\scriptscriptstyle\dagger}}_{\scriptscriptstyle k-1}v^{\scriptscriptstyle\dagger}_{\scriptscriptstyle k-1}
\end{equation}
By (\ref{eq:proofas2}), the second term in (\ref{eq:clt32}) converges to zero almost surely, as $k\rightarrow \infty$. It also converges to zero in expectation, since $\nabla D(\theta)$ is uniformly bounded for $\theta$ in the compact set $\Theta^*$. For the first term in (\ref{eq:clt32}), since the $x_{\scriptscriptstyle k}$ are i.i.d. with distribution $P_{\scriptscriptstyle \theta^*\,}$, it follows that
\begin{equation} \label{eq:clt33}
\mathbb{E}\left[\,u^{\phantom{\scriptscriptstyle{\dagger}}}_{\scriptscriptstyle k}u^{\scriptscriptstyle{\dagger}}_{\scriptscriptstyle k}\,\middle|\mathcal{X}_{\scriptscriptstyle k-1}\right] \,=\, E_{\scriptscriptstyle \theta^*}\left[u(\theta_{\scriptscriptstyle k-1},x)u^{\scriptscriptstyle \dagger}(\theta_{\scriptscriptstyle k-1},x)\right]
\end{equation}
Since $u(\theta,x)$ is a continuous vector field on $\Theta$ for each $x \in X$, and $\theta_{\scriptscriptstyle k-1}$ converge to $\theta^*$ almost surely, it follows that $u(\theta_{\scriptscriptstyle k-1},x)$ converge to $u(\theta^*,x)$ for each $x \in X$, almost surely. On the other hand, it follows from (u2) that the functions under the expectation $E_{\scriptscriptstyle \theta^*}$ in (\ref{eq:clt33}) have bounded second order moments, so they are uniformly integrable~\cite{shiryayev}. Therefore,
\begin{equation} \label{eq:clt34}
\lim\,E_{\scriptscriptstyle \theta^*}\left[u(\theta_{\scriptscriptstyle k-1},x)u^{\scriptscriptstyle \dagger}(\theta_{\scriptscriptstyle k-1},x)\right]\,=\,
E_{\scriptscriptstyle \theta^*}\left[u(\theta^*,x)u^{\scriptscriptstyle \dagger}(\theta^*,x)\right]\,=\,\Sigma^*
\end{equation}
almost surely, by the definition (\ref{eq:cov}) of $\Sigma^*$. It now follows from (\ref{eq:clt32}), (\ref{eq:clt33}), and (\ref{eq:clt34}) that the following limit holds
\begin{equation} \label{eq:clt35}
\lim\,\Sigma_{\scriptscriptstyle k} \,=\, \Sigma^* \hspace{1cm} \text{almost surely}
\end{equation}
To obtain (\ref{eq:clt31}) it is enough to note, as already stated in the verification of Condition (\ref{eq:clt2}), that the $\Sigma_{\scriptscriptstyle k}$ are uniformly bounded in the Frobenius matrix norm. Thus, (\ref{eq:clt35}) implies (\ref{eq:clt31}), by the dominated convergence theorem.  \hfill$\blacksquare$
\end{subequations}

\section{Background on the information metric} \label{sec:efficiency} let $D(\theta)$ be the Kullback-Leibler divergence (\ref{eq:kl}), or any other so-called $\alpha$-divergence~\cite{amari}. Assume the Riemannian metric $\langle\cdot,\cdot\rangle$ of $\Theta$ coincides with the information metric of the model $P$. Then, for any local coordinates $(\tau^{\scriptscriptstyle\alpha}\,;\alpha = 1,\ldots, d\,)$, with origin at $\theta^*$, the following relation holds, by definition of the information metric (see~\cite{amari}, Page 54),
\begin{equation} \label{eq:raofish2}
\left.\frac{\partial^{\scriptscriptstyle\, 2}\!\,D}{\mathstrut\partial\tau^{\scriptscriptstyle \alpha}\partial\tau^{\scriptscriptstyle \beta}}\right|_{\scriptscriptstyle \tau^{\scriptscriptstyle \alpha} = 0} \,=\, \left\langle \frac{\partial}{\partial\tau^{\scriptscriptstyle \alpha}},\frac{\partial}{\partial\tau^{\scriptscriptstyle \beta}}\!\right\rangle_{\scriptscriptstyle\! \theta^*} \\[0.1cm]
\end{equation}
where $\frac{\partial}{\mathstrut\partial\tau^{\scriptscriptstyle\alpha}}$ denote the coordinate vector fields of the local coordinates $\tau^{\scriptscriptstyle\alpha}$. It is also possible to express (\ref{eq:raofish2}) in terms of the Riemannian distance $d(\cdot,\cdot)$, induced by the information metric $\langle\cdot,\cdot\rangle$. Precisely, 
\begin{equation} \label{eq:raofish1}   
  D(\theta) \,=\, \frac{1}{2}\,d^{\scriptscriptstyle\,2}(\theta,\theta^*) \,+\,o\left(d^{\scriptscriptstyle\,2}(\theta,\theta^*)\right)
\end{equation}
This follows immediately from the second-order Taylor expansion of $D(\theta)$, since $\theta^*$ is a minimum of $D(\theta)$, by using (\ref{eq:raofish2}). Formula (\ref{eq:raofish1}) shows that the divergence $D(\theta)$ is equivalent to half the squared Riemannian distance $d^{\scriptscriptstyle\,2}(\theta,\theta^*)$, at $\theta=\theta^*$. 

The scalar products appearing in (\ref{eq:raofish2}) form the components of the information matrix $I^{\scriptscriptstyle\tau}$ of the coordinates $\tau^\alpha$,
$$
I^{\scriptscriptstyle\tau}_{\scriptscriptstyle\alpha\beta} \,=\, 
\left.\frac{\partial^{\scriptscriptstyle\, 2}\!\,D}{\mathstrut\partial\tau^{\scriptscriptstyle \alpha}\partial\tau^{\scriptscriptstyle \beta}}\right|_{\scriptscriptstyle \tau^{\scriptscriptstyle \alpha} = 0}
$$
In any change of coordinates, these transform like the components of a $(0,2)$ covariant tensor~\cite{petersen}. That is, if $(\theta^{\scriptscriptstyle\,\alpha}\,;\alpha = 1,\ldots, d\,)$ are any local coordinates defined at $\theta^*$,
$$
I^{\scriptscriptstyle\tau}_{\scriptscriptstyle\alpha\beta} \,=\, 
\left(\frac{\partial \theta^{\scriptscriptstyle\,\gamma}}{\partial\tau^{\scriptscriptstyle\alpha}}\right)_{\!\scriptscriptstyle \theta^*}
I^{\scriptscriptstyle\theta}_{\scriptscriptstyle\gamma\kappa}
\left(\frac{\partial \theta^{\scriptscriptstyle\,\kappa}}{\partial\tau^{\scriptscriptstyle\beta}}\right)_{\!\scriptscriptstyle \theta^*}
$$ 
where the subscript $\theta^*$ indicates the derivative is evaluated at $\theta^*$, and where $I^{\scriptscriptstyle\theta}_{\scriptscriptstyle\gamma\kappa}$ are the components of the information matrix $I^{\scriptscriptstyle \theta}$ of the coordinates $\theta^{\scriptscriptstyle\,\alpha}$. 

The recursive estimates $\theta_{\scriptscriptstyle n}$ are said to be asymptotically efficient, if they are asymptotically efficient in any local coordinates $\tau^{\scriptscriptstyle\alpha}$, with origin at $\theta^*$. That is, according to the classical definition of asymptotic efficiency~\cite{ibrahas}\cite{vaart}, if the following weak limit of probability distributions is verified~\cite{shiryayev},
\begin{equation} \label{eq:efficiency}
  \mathcal{L} \left\lbrace (n^{\scriptscriptstyle 1/2}\tau^{\scriptscriptstyle\alpha}_{\scriptscriptstyle n})\right\rbrace \,\Longrightarrow\, N_{\scriptscriptstyle d}\left(0,\Sigma^{\scriptscriptstyle\tau}\right) \hspace{1cm} \Sigma^{\scriptscriptstyle\tau}\,=\,\left(I^{\scriptscriptstyle \tau}\right)^{\scriptscriptstyle -1}
\end{equation}
where $\mathcal{L}\lbrace\ldots\rbrace$ denotes the probability distribution of the quantity in braces, $\tau^{\scriptscriptstyle \alpha}_{\scriptscriptstyle n} = \tau^{\scriptscriptstyle \alpha}(\theta_{\scriptscriptstyle n})$ are the coordinates of the recursive estimates $\theta_{\scriptscriptstyle n\,}$, and $N_{\scriptscriptstyle d}\left(0,\Sigma^{\scriptscriptstyle\tau}\right)$ denotes a centred $d$-variate normal distribution with covariance matrix $\Sigma^{\scriptscriptstyle\tau}$.  

It is important to note that asymptotic efficiency of the recursive estimates $\theta_{\scriptscriptstyle n}$ is an intrinsic geometric property, which does not depend on the particular choice of local coordinates $\tau^{\scriptscriptstyle\alpha}$, with origin at $\theta^*$. This can be seen from the transformation rule of the components of the information matrix, described above. In fact, since these transform like the components of a $(0,2)$ covariant tensor, the components of $\Sigma^{\scriptscriptstyle\tau}$ transform like those of a $(2,0)$ contravariant tensor, which is the correct transformation rule for the components of a covariance matrix. 


\bibliographystyle{siamplain}
\bibliography{references}

\end{document}